\documentclass[preprint,11pt]{elsarticle}

\usepackage[singlelinecheck=false,labelsep=period,font=footnotesize]{caption}
\usepackage[scale=0.66]{geometry}

\usepackage{comment}
\usepackage{pdfpages}
\usepackage{mathtools}
\usepackage{amsmath, amsthm, amssymb, amsfonts}
\usepackage{empheq}  
\usepackage{lineno,hyperref}
\usepackage{comment}
\usepackage{multirow, array}
\usepackage{graphicx}
\usepackage{float}
\usepackage{color}
\usepackage{multicol}

\usepackage{multirow}

\usepackage{xr-hyper}

\usepackage{hyperref}

\usepackage{graphicx}

\usepackage{amsfonts}

\usepackage{amsmath, amsthm}

\usepackage{color}

\usepackage[normalem]{ulem} 

\usepackage{epstopdf}

\usepackage{subcaption}

\usepackage{eurosym}

\usepackage{enumitem}

\usepackage{bbm}

\usepackage{pgfplots,pgfplotstable,booktabs}

\pgfplotsset{compat=1.8}

\usepgfplotslibrary{statistics}

\usepackage{balance}

\allowdisplaybreaks

\newtheorem{thm}{Theorem}

\newtheorem{remark}{Remark}
\newtheorem{example}{Example}
\newtheorem{prop}[thm]{Proposition}

\newtheorem{definition}{Definition}

\newenvironment{ldescription}[1]
  {\begin{list}{}%
   {\renewcommand\makelabel[1]{##1\hfill}%
   \settowidth\labelwidth{\makelabel{#1}}%
   \setlength\leftmargin{\labelwidth}
   \addtolength\leftmargin{\labelsep}}}
  {\end{list}}

\usepackage{lineno,hyperref}

\makeatletter
\def\ps@pprintTitle{%
  \let\@oddhead\@empty
  \let\@evenhead\@empty
  \let\@oddfoot\@empty
  \let\@evenfoot\@oddfoot
}
\makeatother
\begin{document}

\begin{frontmatter}

\title{Distributionally Robust Optimal Power Flow with Contextual Information}

\author{Adri\'an Esteban-P\'erez, Juan M. Morales\corref{mycorrespondingauthor}}
\address{Department of Applied Mathematics,  University of M\'alaga, M\'alaga, 29071, Spain}


\cortext[mycorrespondingauthor]{Corresponding author}
\ead{adrianesteban@uma.es, juan.morales@uma.es}


\begin{abstract}
In this paper, we develop a distributionally robust chance-constrained formulation of the Optimal Power Flow problem (OPF)  whereby the system operator can leverage contextual information. For this purpose, we exploit an ambiguity set based on probability trimmings and optimal transport through which the dispatch solution is protected against the incomplete knowledge of the relationship between the OPF uncertainties and the context that is conveyed by a sample of their joint probability distribution. We provide {\color{black} a tractable} reformulation of the proposed distributionally robust chance-constrained OPF problem under the popular conditional-value-at-risk approximation. By way of numerical experiments run on a modified IEEE-118 bus network with wind uncertainty, we show how the power system can substantially benefit from taking into account the well-known statistical dependence between the point forecast of wind power outputs and its associated prediction error. Furthermore, the  experiments conducted  also reveal that the distributional robustness conferred on the OPF solution by our probability-trimmings-based approach is superior to that bestowed by alternative approaches in terms of expected cost and system reliability.
\end{abstract}


\end{frontmatter}


\section{Introduction}
\label{sec:introduction}

The Optimal Power Flow (OPF) is a fundamental problem in power system operations. Traditionally, the goal of the OPF problem is to minimize the cost of the power generation dispatch that supplies the electricity demand while complying with some physical and engineering constraints. The growing penetration of electricity generation sources like wind and solar power, which are intrinsically uncertain, has led power system engineers to account for randomness in OPF analyses. Hence, the OPF is to be formulated today as an optimization problem \emph{under uncertainty}.

A common way to cope with uncertainty in the constraints of an optimization problem and, in particular, of an OPF model is by way of the so-called \emph{chance constraints}, which allow the modeler to impose the constraint satisfaction with a certain probability only. Accordingly, chance-constrained optimal power flow models (CC-OPF) have been developed to control the violation probability of, for instance, line and generation capacity limits. In particular, references \cite{jia2021iterative, pena2020dc, vrakopoulou2013probabilistic} consider \emph{joint} chance constraints, by which the system operator enforces that all constraints must simultaneously hold with a probability greater than or equal to $1-\epsilon$, where $\epsilon \in (0,1)$ is a pre-fixed acceptable tolerance of dispatch infeasibility. This is in contrast to \emph{single} chance constraints, whereby the probability of constraint satisfaction is imposed on each constraint of the OPF separately. Although single chance constraints have been considered in the  literature of CC-OPF models due to their attractive tractability properties (see, e.g., \cite{Bienstock2014, Lubin2016,zhang2011chance} and references therein), {\color{black} they do not provide guarantees on the \emph{joint} satisfaction of the OPF constraints in general or, when they do, they usually lead to costly and over-conservative dispatch solutions \cite{pena2020dc}}.

One of the main challenges in solving CC-OPF problems is that the underlying probability distribution of the random variables affecting the OPF constraints is generally unknown. In fact, in practice, only past historical observations of those variables are available to the system operator. Within this context, Distributionally Robust Optimization (DRO) has emerged as a versatile and powerful paradigm to tackle the ambiguity of the distribution of an optimization problem's random parameters.  DRO, {\color{black}therefore} considers a set of potential distributions for those parameters, namely, the so-called \emph{ambiguity set} \cite{Rahimian2019}. In the literature on DRO, there exist two distinct ways to specify an ambiguity set, either using statistical moments or probability metrics.
%
%
The main drawback of using a \emph{moment-based} DRO approach is that the available historical information is only used to estimate the moments of the true data-generating distribution and thus, if new samples of said distribution become available, but the moment estimates are not modified, the ambiguity set does not change. Alternatively, in \emph{metric-based} DRO, the ambiguity set collects all potential distributions whose distance to a particular nominal distribution is lower than or equal to some pre-fixed value. Although there are different choices of probability metrics to construct the ambiguity set, the Wasserstein distance has received much attention within the power systems community. This is partly due to the performance guarantees provided by this distance (see, e.g., \cite{MohajerinEsfahani2018}) in contrast to the moment-based DRO approach.

Thus, distributionally robust chance-constrained optimal power flow  (DRCC-OPF) models seek the optimal dispatch such that all the model constraints are satisfied with a pre-fixed confidence level for all the probability distributions within the ambiguity set built either by means of the Wasserstein metric  (\cite{arrigo2021,arrigo2022,Ordoudis2021,Zhou2020}), using moments (\cite{Li2016dro, Li2019,Li2019ieeecontrol,Xie2018ieee}), or by way of a  discrete probability distribution with probability masses and locations varying within a box~\cite{Jabr2020}.

 Unfortunately, the consideration of \emph{joint} chance constraints in a DRO framework (\cite{Hanasusanto2017OPRECHANCE,Xie2020OPRECHANCE}) renders an intractable optimization problem in general. For this reason, researchers have considered DRO-OPF models based on the well-known conservative, but far more tractable approximation, given by the concept of \emph{conditional-value-at-risk} (${\mathbf{CVaR}}$), see, for example, \cite{arrigo2021,Jabr2020,Ordoudis2021} and references therein. The authors in \cite{Chen2018b} show that the  ${\mathbf{CVaR}}$ offers a tight convex approximation of the chance constraints under the DRO framework, which justifies its popularity. In addition, by way of ${\mathbf{CVaR}}$-based chance-constraints, the power system operator can control not only the violation  probability, but also the violation magnitude, which can be important from the standpoint of power system operations.

Yet another key issue is that the system operator is able to not only  tune the robustness of the resulting chance-constrained DRO model by adjusting the probability of constraint violation, but also, and very importantly, through the specificity degree of the ambiguity set. In this regard, some previous approaches have focused on producing more meaningful ambiguity sets by incorporating structural information on the underlying true probability distribution (see, e.g., \cite{arrigo2021, arrigo2022, Li2016dro,Li2019,Li2019ieeecontrol}). {\color{black} More specifically, the authors in \cite{arrigo2021} introduce a DRCC-OPF   model with single chance constraints, which considers all distributions within a Wasserstein ball that conform to a given copula-based dependence structure among the random variables. The authors in \cite{arrigo2022} propose
an iterative algorithm for solving a bilinear exact reformulation of a DRCC-OPF model with single chance constraints considering some support information. The authors of
\cite{Li2016dro} and \cite{Li2019} provide and study DRCC-OPF models where some moment and unimodality information is included in the ambiguity set. Both aforementioned approaches are extended in~\cite{Li2019ieeecontrol} to allow for misspecified modes.
}

{\color{black}
Ideally, one would like to have the smallest ambiguity set that contains the true data-generating distribution. In this vein, if we have some information on the true distribution, we should use it to discard all those other distributions that do not conform with that information from the ambiguity set. As mentioned above, that information can be, for example, some dependence structure via copulas  \cite{arrigo2021}, support information \cite{arrigo2022} or shape  information (such as unimodality) \cite{Li2019,Li2019ieeecontrol}. Our aim is also to leverage a more informed ambiguity set,  but, and for the first time to our knowledge, that information refers to a given \emph{context}.
This side/contextual information (also known as \emph{features} or \emph{attributes}) is related to outcomes of random variables that may have predictive power on the OPF's uncertainties. Accordingly, we make use of an ambiguity set that accounts for the possible dependence between the uncertainties and these explanatory variables.
Thus, the contextual information allows us to discard implausible distributions.
}

 More specifically, in the work we present here, we  exploit the contextual information provided by the \emph{point forecasts} of those uncertainties. In the energy forecasting community,it is a well known fact that the power forecast error of a wind farm highly depends on the wind power forecast itself \citep{Bludszuweit2008,Fabbri2005}. Within the context of DRCC-OPF, this means that the wind power point forecast constitutes valuable information to build a proper ambiguity set for the wind power forecast error.


The main contributions of this work are thus:
\begin{enumerate}

\item We provide a formulation of the DC-OPF problem with \emph{joint} chance constraints as a \emph{conditional stochastic program}, in which both the expected cost of the power dispatch and the chance constraints are conditioned on some contextual information, in particular, the point forecasts of the OPF's uncertainties. {\color{black}To the best of our knowledge, this work is the first to tackle a chance constraint system with a distributionally robust approach \emph{that accounts for contextual information}.}

\item To tackle the resulting conditional stochastic program, we propose an entirely data-driven and non-parametric DRO approach, in which no assumption on the relationship between the context and the OPF's uncertainties (i.e., the point forecasts and their error in our case) is made. {\color{black} Our approach makes use of} the ambiguity set based on probability trimmings that is introduced in \cite{Estebanmorales_trimmings_extended2020} and which is able to integrate contextual information. We prove that, under the widely used ${\mathbf{CVaR}}$-based approximation of the joint chance constraint system \cite{nemirovski2007convex}, it results in a tractable distributionally robust joint chance-constrained OPF model.

\item Finally, by way of numerical experiments in which we compare our DRCC-OPF formulation with {\color{black} alternative approaches in the literature}, we show that exploiting the contextual information provided by the wind power point forecasts allows identifying dispatch solutions with a better trade-off between expected cost and system reliability.
\end{enumerate}

The remainder of this paper is organized as follows. Section~\ref{sec:model_formulation} introduces the mathematical formulation of the distributionally robust joint chance-constrained DC-OPF problem with contextual information that we propose. Then, Section~\ref{sec:CVaRapprox} provides a tractable reformulation of the distributionally robust joint chance constraints under the well-known ${\mathbf{CVaR}}$ approximation, while Section~\ref{sec:worstEC} focuses on reformulating the worst-case expected cost in a manageable way. Results from numerical experiments are presented and discussed in Section~\ref{sec:numerics}. Finally, Section~\ref{sec:conclusion} concludes the paper with some final remarks.  The manuscript  also contains appendices with the main notation used throughout the main text and supportive optimization models. {\color{black}The proofs of theoretical results and additional numerical experiments are available in the supplementary material associated with this article.}\\
{\color{black}\emph{Notation}.
Unless otherwise stated, in this paper, we use boldface lowercase letters to represent arrays and boldface capital letters for matrices. Vector $\mathbf{1}$ ($\mathbf{0}$) is a vector of ones (zeros) of appropriate dimension,  and the inner product of two vectors $\mathbf{u}, \mathbf{v}$ will be denoted by $\langle \mathbf{u}, \mathbf{v}\rangle$. The cardinality of a set $A$ will be indicated by $|A|$.  The support function of a set $B \subseteq \mathbb{R}^d$, $S_B$, is defined as $S_B(a):=\sup_{b \in B }\langle a,b\rangle$. Moreover, we reserve the symbol ``$\;\widehat{\,}\;$" for objects which are dependent on the sample data and denote ``expectation" with the symbol $\mathbb{E}$. In addition, throughout the paper we assume that we always have measurability for those objects whose expected values we consider. The reader is referred to \ref{appendix_notation} for a complete list of the notation used throughout this paper.

}
\section{DC-OPF under uncertainty: Mathematical Formulation}\label{sec:model_formulation}

{\color{black}
Next we introduce the DC-OPF problem under uncertainty. The problem is formulated as a distributionally robust version of the \emph{joint} chance-constrained DC-OPF model described in~\cite{pena2020dc}, where we have also accounted for the procurement of reserve capacity and its associated cost, as in \cite{Li2019, Li2019ieeecontrol}. Nonetheless, unlike in~\cite{pena2020dc}, where the generators' cost functions are assumed to be quadratic, here we model those costs as convex piecewise linear functions. Furthermore, there exists a number of different variants of the distributionally robust chance-constrained DC-OPF problem (e.g., \cite{arrigo2022,Jabr2020,Li2019,Ordoudis2021,Zhou2020}), which essentially differ in the treatment of the chance constraints (single vs. joint), the cost structure of generators that is assumed, and the ambiguity set used. What makes our formulation unique among those variants is its ability to exploit contextual information.

\subsection{Variables and constraints}

Consider a power system with a set $\mathcal{L}$ of transmission lines, a set $\mathcal{B}$ of buses, a set  $\mathcal{W}$ of wind power plants (or, more generally, weather-dependent renewable generators), and  a set  $\mathcal{G}$ of conventional generators (i.e., dispatchable units that are not weather-dependent). For ease of formulation, power loads are assumed to be deterministic. Next we introduce  each of the main components of the DC-OPF problem.

\begin{enumerate}
    \item \emph{Wind power plants.} For each wind power plant $m \in \mathcal{W}$, the random power output is modeled as $f_m+\omega_m$, where $f_m$ is the predicted power output and $\omega_m$ is the (random) wind forecast error at wind power plant $m$. We denote  the system-wise aggregate wind power forecast error as $\Omega$, i.e.,   $\Omega:=\sum_{m\in \mathcal{W}}\omega_m=\langle \mathbf{1},\boldsymbol{\omega} \rangle $. Let
    $\mathbf{f}:=(f_m)_{m\in \mathcal{W}}, \boldsymbol{\omega}:=(\omega_m)_{m\in \mathcal{W}} $ be the array of predicted power outputs and wind power prediction errors, respectively.
    \item \emph{Generators:}  For each $j\in \mathcal{G}$,  the actual power output of generator $j$, $\tilde{g}_j(\boldsymbol{\omega})$,  is expressed as the sum of the scheduled generation, $g_j$, and the (random) adjusted  power $\tilde{r}_j(\boldsymbol{\omega})$ (also known as \emph{deployed reserve}).
   As customary, we assume an affine control policy to counterbalance the  wind forecast errors by deploying generators' reserves \cite{Bienstock2014}, that is,
    \begin{equation}\label{reserv_affine}
           \tilde{g}_j(\boldsymbol{\omega}):=g_j+\tilde{r}_j(\boldsymbol{\omega})=g_j-\beta_j\Omega=g_j-\beta_j\langle  \mathbf{1}, \boldsymbol{\omega}  \rangle,\enskip \forall j \in \mathcal{G}
    \end{equation}
    where $\beta_j$ is the participation factor of generator $j$. Denote by $\boldsymbol{\beta}:=(\beta_j)_{j\in \mathcal{G}}, \mathbf{g}:=(g_j)_{j\in \mathcal{G}}$ the array of non-negative participation factors and scheduled generation, respectively. Let $\tilde{\mathbf{g}}(\boldsymbol{\omega}):=(\tilde{g}_j(\boldsymbol{\omega}))_{j \in \mathcal{G}}$, $ \tilde{\mathbf{r}}(\boldsymbol{\omega}):=(\tilde{r}_j(\boldsymbol{\omega}))_{j \in \mathcal{G}}=(-\beta_j\langle  \mathbf{1}, \boldsymbol{\omega}  \rangle)_{j \in \mathcal{G}}$  be the array of actual power outputs and deployed reserves, in that order.

   The  following constraints determine the provision of reserve capacities:
    \begin{equation}\label{reserve_limits}
         -\mathbf{r}^D\leqslant\tilde{\mathbf{r}}(\boldsymbol{\omega}) \leqslant \mathbf{r}^U
    \end{equation}
    with $ \mathbf{r}^{D}, \mathbf{r}^{U}$ being the arrays of downward and upward reserve capacity provided by the generators, respectively.

 Naturally, the  following technical constraints, which link the generation dispatches and the provision of reserve capacities, must hold:
      \begin{align}
       &  \mathbf{g}+\mathbf{r}^U \leqslant  \mathbf{g}^{\max},  \label{reserve_limits_u_phys} \\
      & \mathbf{g}-\mathbf{r}^D\geqslant   \mathbf{g}^{\min} \label{reserve_limits_d_phys}
    \end{align}
    where $  \mathbf{g}^{\min},  \mathbf{g}^{\max}$ are the arrays of minimum and maximum power output of the generators, respectively.

    \item \emph{Network constraints.} The total power generation must equal the total system demand (\emph{power balance constraint}), that is,
    \begin{equation}\label{balance_random}
    \langle \mathbf{1},\tilde{\mathbf{g}}(\boldsymbol{\omega})  \rangle +\langle \mathbf{1}, \mathbf{f}+\boldsymbol{\omega} \rangle
       =\langle \mathbf{1}, \mathbf{L} \rangle
    \end{equation}
    where  $\mathbf{L}:=(L_b)_{b\in \mathcal{B}}$ denotes the array of nodal loads.  Using \eqref{reserv_affine}, Eq.~\eqref{balance_random} is equivalent to:
    \begin{align}
       \langle \mathbf{1},\mathbf{g}  \rangle +\langle \mathbf{1}, \mathbf{f} \rangle&=\langle \mathbf{1}, \mathbf{L}   \rangle \label{balance_cons} \\
        \langle \mathbf{1}, \boldsymbol{\beta} \rangle&=1, \enskip \boldsymbol{\beta} \geqslant \mathbf{0}\label{particip_cons}
    \end{align}
which guarantee the power balance both in the dispatch and the real-time stages, respectively.\\
Finally, we assume that the power flow through the lines is given by a linear function of the nodal power injections, that is, $\mathbf{M}^{\mathcal{G}}(\tilde{\mathbf{g}}(\boldsymbol{\omega}))+\mathbf{M}^{\mathcal{W}}(\mathbf{f}+\boldsymbol{\omega})-\mathbf{M}^{\mathcal{B}}\mathbf{L}$, based on the DC power flow approximation, where  $\mathbf{M}^{\mathcal{G}}, \mathbf{M}^{\mathcal{W}}$ and $\mathbf{M}^{\mathcal{B}}$ are the matrix for generators, wind plants and loads  given by  the DC power transfer distribution factors   \cite{Stott2009}, in that order. Hence,
 the constraints
  \begin{align}\label{flow_limitd}
        -\textbf{Cap}  \leqslant & \mathbf{M}^{\mathcal{G}}(\tilde{\mathbf{g}}(\boldsymbol{\omega}))+\mathbf{M}^{\mathcal{W}}(\mathbf{f}+\boldsymbol{\omega})-\mathbf{M}^{\mathcal{B}}\mathbf{L} \leqslant \textbf{Cap}
  \end{align}
  enforce the transmission capacity limits  where
  $\textbf{Cap} :=(\text{Cap}_{\ell})_{\ell \in \mathcal{L}} $ denotes the array of transmission line capacities.
  \end{enumerate}
  }

  \subsection{Dealing with uncertainty in the DC-OPF problem}\label{sec:OPFuUnc}
%

%
In practice, it is often the case that the random vector of forecast errors $\boldsymbol{\omega}$ shows some statistical dependence on some features/covariates, which we can model, in general, by some random vector $\mathbf{z}$. In fact, the forecast wind power output $\mathbf{f}$ serves in itself as an obvious explanatory random vector for the subsequent forecast error $\boldsymbol{\omega}$. In this approach, we want to exploit this statistical dependence to identify a better power generation dispatch and provision of reserve capacity.

Let $\mathbf{z}:={\color{black}(z_m)_{m\in \mathcal{W}}}$ be the random vector modeling the features and let
$\mathbb{Q}$ be the probability measure of the joint distribution of $\mathbf{z}$ and $\boldsymbol{\omega}$, which is supported on $\Xi$. {\color{black} For convenience, we define $\boldsymbol{\xi}:=(\mathbf{z}, \boldsymbol{\omega})$}.   Given the array of forecast wind power outputs,  $\mathbf{f}:={\color{black}(f_m)_{m\in \mathcal{W}}}$, {\color{black}set the contextual information $\boldsymbol{\xi}:=(\mathbf{z},\boldsymbol{\omega})\in \widetilde{\Xi}$}  defined by  the event $( {\color{black}\mathbf{z}=\mathbf{f}};\  \boldsymbol{\omega}\in \widetilde{\Xi}_{\boldsymbol{\omega}})$, with $\widetilde{\Xi}_{\boldsymbol{\omega}}$ being the support of $\boldsymbol{\omega}$ conditional on ${\color{black}\mathbf{z}=\mathbf{f}}$. {\color{black} The errors of forecasting the power output of a wind farm are naturally bounded. Their lower bound is the forecast value itself, while their upper bound is given by the difference of the capacity of the wind farm and the predicted value. Therefore,} $\widetilde{\Xi}_{\boldsymbol{\omega}}$ is the hypercube $\prod_{m {\color{black}\in \mathcal{W}}}  [-f_m,\overline{C}_m-f_m]$, where $\overline{C}_m$ represents the capacity of wind farm $m$. {\color{black}Note that the optimal   dispatch  is, therefore, parametrized on the predicted wind power outputs $\{f_m \}_{m\in\mathcal{W}}$. }




In real life, however, neither the joint distribution $\mathbb{Q}$, nor the conditional one $\mathbb{Q}_{\boldsymbol{\omega}/\mathbf{z}=\mathbf{f}}$, are known. The system operator only has  access to a finite set of samples of size $N$ (i.e. the training set) of the true joint distribution $\mathbb{Q}$, which we denote as
{\color{black}$\widehat{\Xi}^N_{\omega}:=\{\widehat{\boldsymbol{\xi}_i} \}_{i=1}^N=\{(\widehat{\mathbf{z}}_i,\widehat{\boldsymbol{\omega}}_i)\}_{i=1}^N$}. In our context, $\widehat{\Xi}^N_{\omega}$ is made up of $N$ past observations of the predicted wind power outputs and their associated errors. Hence, the system operator needs to infer or construct a proxy of $\mathbb{Q}_{\boldsymbol{\omega}/\mathbf{z}=\mathbf{f}}$ from the sample $\widehat{\Xi}^N_{\omega}$, so that this proxy can be used to compute a reliable and cost-efficient OPF solution. However, the limited information that $\widehat{\Xi}^N_{\omega}$ conveys on $\mathbb{Q}_{\boldsymbol{\omega}/\mathbf{z}=\mathbf{f}}$ makes this inference process ambiguous, and as such, we propose  employing the following \emph{distributionally robust} chance-constrained OPF model to protect the system operator's decision against this ambiguity:
{\color{black}
 \begin{align}
\min_{\mathbf{x}\in X}\; &  \sup_{Q_{\widetilde{\Xi}} \in \widehat{\mathcal{U}}} \mathbb{E}_{Q_{\widetilde{\Xi}}} \left[  C(\tilde{\mathbf{g}}(\boldsymbol{\omega}))+\langle \mathbf{c}^D, \mathbf{r}^D \rangle +\langle \mathbf{c}^U, \mathbf{r}^U \rangle \right] \label{worst_case_expected} \\
\text{s.t.}\;  &\inf_{Q_{\widetilde{\Xi}} \in \widehat{\mathcal{U}}} Q_{\widetilde{\Xi}}
\left( \begin{array}{c}
         -\mathbf{r}^D\leqslant \tilde{\mathbf{r}}(\boldsymbol{\omega}) \leqslant \mathbf{r}^U\\
         -\textbf{Cap}  \leqslant  \mathbf{M}^{\mathcal{G}}(\tilde{\mathbf{g}}(\boldsymbol{\omega}))+\mathbf{M}^{\mathcal{W}}(\mathbf{f}+\boldsymbol{\omega})-\mathbf{M}^{\mathcal{B}}\mathbf{L} \leqslant \textbf{Cap}
 \end{array} \right) \geqslant 1-\epsilon \label{chance_constraintsDRO}
\end{align}
 where $X$ is  the deterministic feasible set for the array of decision variables $ \mathbf{x}=(\mathbf{g},\boldsymbol{\beta},\mathbf{r}^D,\mathbf{r}^U)$ defined by the constraints        \eqref{reserve_limits_u_phys},  \eqref{reserve_limits_d_phys}, \eqref{balance_cons} and  \eqref{particip_cons}.

 {\color{black} The set
 $\widehat{\mathcal{U}}$  in \eqref{worst_case_expected}-\eqref{chance_constraintsDRO}} stands for an ambiguity set for the true conditional distribution $\mathbb{Q}_{\boldsymbol{\omega}/\mathbf{z}=\mathbf{f}}$. Here we will use the ambiguity set based on probability trimmings and optimal transport introduced in   \cite{Estebanmorales_trimmings_extended2020}, which allows us  to exploit the side information within a DRO framework  in a fully data-driven sense. This ambiguity set is formally introduced in Section~\ref{sec:amb_set}.

Objective function~\eqref{worst_case_expected}  minimizes the    expected total operational cost over the worst-case distribution from the ambiguity set $\widehat{\mathcal{U}}$, including  the sum of the (random) total generation cost,  $
        C(\tilde{\mathbf{g}}(\boldsymbol{\omega}))
    $,   and the (deterministic) cost  of providing  up- and  down-reserve  capacities, $\langle \mathbf{c}^D, \mathbf{r}^D \rangle +\langle \mathbf{c}^U, \mathbf{r}^U \rangle$. The total generation cost function $C(\cdot)$ is assumed to be given
    by the sum of $|\mathcal{G}|$ convex piecewise linear cost functions  with $S_j$ pieces/blocks, i.e.,
    $C(\tilde{\mathbf{g}}(\boldsymbol{\omega})):=\sum_{j \in \mathcal{G}} \max_{s=1,\ldots, S_j } \{ m_{js}\tilde{g}_j(\boldsymbol{\omega}) +n_{js} \}$,
    where $m_{js}, n_{js}$ stand for  the slope and the intercept of the $s$-th piece   for generator $j$, respectively. Parameters
  $\mathbf{c}^D, \mathbf{c}^U$ are the arrays of  downward and upward reserve procurement cost of the generators, respectively. Note that the wind power production cost is assumed to be zero.
Finally, constraint \eqref{chance_constraintsDRO}
 establishes a tolerance $\epsilon$  in terms of the joint violation probability of the OPF constraints under any conditional probability distribution $Q_{\widetilde{\Xi}}$ for $\boldsymbol{\omega}$ (given some contextual information) within  the ambiguity set $ \widehat{\mathcal{U}}$, which will be defined right after the following remark.

\begin{remark}
Estimating the probability distribution of the wind power forecast error conditional on the point prediction, i.e., $\mathbb{Q}_{\boldsymbol{\omega}/\mathbf{z}=\mathbf{f}}$  (a task that is typically referred to as \emph{probabilistic forecasting}) requires \emph{postulating a statistical model for that distribution}, following the scheme
{\color{black} ``from data to decision through prediction''  }.
The ``prediction'' step in that scheme may negatively affect the decision value in the end due to model misspecification, the estimation error and/or simply because of the so-called  ``optimizer’s curse” or ``optimization bias” \cite{MohajerinEsfahani2018}.

{\color{black} Our approach, in contrast, is fully data driven in the sense that it follows an
alternative scheme by which the decision is directly inferred from the data without being forced to use  a specific predictive model (which could be the wrong choice).} Furthermore, as we explain below, our approach protects the decision against the uncertainty intrinsic to the process of inferring the \emph{conditional} forecast error distribution  from samples of the \emph{joint} distribution of the pair (point prediction, prediction error). This can be particularly important when the sample size is small and said uncertainty is high as a result.
\end{remark}

}

 {\color{black}

\subsection{The ambiguity set based on probability trimmings} \label{sec:amb_set}
To formalize the ambiguity set based on probability trimmings that will play the role of $\widehat{\mathcal{U}}$ in \eqref{worst_case_expected}-\eqref{chance_constraintsDRO}, we first need to introduce a series of probability-related concepts. We start with the Wasserstein metric of order $1$, a well-known  metric for probability distributions with finite first moment closely linked to the optimal transport problem~\cite{Santambrogio2015}:

 \begin{definition}[$1$-Wasserstein metric]\label{Wass1}
 Given two probability measures $P,Q$ with finite first moment supported on $\Xi$, that is,   $P,Q \in \mathcal{P}_1(\Xi)$, the  Wasserstein metric of order~$1$ between $P$ and $Q$, $W_1\big (P,Q\big ) $, is given by the value
$$\begin{aligned} \inf_{\Pi} \left\{ \int _{\Xi^2} \Vert \boldsymbol{\xi}_1 - \boldsymbol{\xi} _2 \Vert_1 \, \Pi (\mathrm {d}\boldsymbol{\xi} _1, \mathrm {d}\boldsymbol{\xi} _2) ~: \begin{array}{l} \Pi \textit{ is a joint distribution of } \boldsymbol{\xi}_1 \textit{ and } \boldsymbol{\xi}_2 \\ \textit{with marginals } P \textit{ and } Q, \textit{ respectively} \end{array}\right\} \end{aligned}$$
 \end{definition}

Next, we introduce the definition of trimmings and trimming sets of an empirical probability measure.  Essentially, trimming a probability measure allows to play down the weight of some regions of the measurable
space without completely removing them from the support~set.

\begin{definition}[$(1-\alpha)$-empirical trimmings sets]\label{def_trimming}
Consider the sample data $\{\widehat{\boldsymbol{\xi}}_i\}_{i=1}^{N}$ and their associated empirical measure $\widehat{\mathbb{Q}}_{N} = \frac{1}{N}\sum_{i=1}^{N}{\delta_{\widehat{\boldsymbol{\xi}}_i}}$.
 If $\alpha > 0$, the set of all $(1-\alpha)$-trimmings of $\widehat{\mathbb{Q}}_{N}$ is given by all probability distributions in the form $\sum_{i=1}^{N}{b_i\delta_{\widehat{\boldsymbol{\xi}}_i}}$ such that $0\leq b_i\leq \frac{1}{N\alpha}$, $\forall i =1, \ldots, N$, and $\sum_{i=1}^{N}{b_i}=1$.
\end{definition}

For ease of understanding, we provide below an example of a $(1-\alpha)$-empirical trimmings set.

\begin{example}
Suppose the empirical joint measure $\widehat{\mathbb{Q}}_N:= \sum_{i=1}^{3}{\delta_{(\widehat{z}_i,\widehat{\omega}_i)}} = \frac{1}{3}(\delta_{(1,0)}+ \delta_{(0,5)}+ \delta_{(2,3)})$ ($N = 3$).
If $\alpha = 0.5$, then $\frac{1}{N\alpha} = \frac{1}{3\cdot 0.5} = \frac{2}{3}$. Therefore, the $0.5$-trimmings set of $\widehat{\mathbb{Q}}_N$ is given by

$$  \mathcal{R}_{0.5}(\widehat{\mathbb{Q}}_N)  :=\Bigg\{\sum_{i=1}^{3}{b_i\delta_{\widehat{\boldsymbol{\xi}}_i}} : 0\leq b_i\leq \frac{2}{3}, \forall i =1, \ldots, 3; \sum_{i=1}^{3}{b_i}=1\Bigg\}$$
The following statements hold thus true:
\begin{align*}
 \widehat{\mathbb{Q}}_N&=\frac{1}{3}(\delta_{(1,0)}+ \delta_{(0,5)}+ \delta_{(2,3)}) \in \mathcal{R}_{0.5}(\widehat{\mathbb{Q}}_N)   \\
 \mathcal{Q}&=\frac{2}{3}\delta_{(1,0)}+ \frac{1}{3}\delta_{(0,5)} \in \mathcal{R}_{0.5}(\widehat{\mathbb{Q}}_N)  \\
  \mathcal{P} &=  \delta_{(1,0)} \not \in \mathcal{R}_{0.5}(\widehat{\mathbb{Q}}_N) \; \textrm{(the trimming must retain one point and a half at least)}\\
 \mathcal{S} &= \frac{2}{3}\delta_{(1,0)}+ \frac{1}{6}\delta_{(0,5)} + \frac{1}{6} \delta_{(2,3)} \in \mathcal{R}_{0.5}(\widehat{\mathbb{Q}}_N)\\
  \mathcal{V}  &=\frac{3}{4}\delta_{(1,0)}+ \frac{1}{12}\delta_{(0,5)} + \frac{2}{12} \delta_{(2,3)} \not\in \mathcal{R}_{0.5}(\widehat{\mathbb{Q}}_N) \;  (because\ b_1 > 2/3)
\end{align*}

\end{example}


Our ambiguity set also relies on the concept of \emph{minimum transportation budget}, which is provided below.


\begin{definition}[Minimum transportation budget]\label{MTB}
Given $\alpha >0$, the \emph{minimum transportation budget}, which we denote as $\underline{\rho}_{N\alpha}$, is the $1$-Wasserstein distance  between the set of probability distributions $\mathcal{P}_1(\widetilde{\Xi})$ and the $(1-\alpha)$-trimming of the empirical distribution $\widehat{\mathbb{Q}}_N$ that is the \emph{closest} to that set, that is,

\begin{equation}\label{def_MTB}
    \underline{\rho}_{N\alpha} = \frac{1}{N\alpha}\sum_{k=1}^{\lfloor N\alpha \rfloor}{\textrm{dist}(\boldsymbol{\xi}_{k:N}, \widetilde{\Xi})} + \left(1-\frac{\lfloor N\alpha \rfloor}{N\alpha}\right) \textrm{dist}(\boldsymbol{\xi}_{\lceil N\alpha\rceil:N}, \widetilde{\Xi})
\end{equation}
where $\boldsymbol{\xi}_{k:N}$ is the $k$-\emph{th} nearest data point from the sample to set $\widetilde{\Xi}$ and
$
\textrm{dist}
    (\boldsymbol{\xi}_j,
    \widetilde{\Xi}):=
    \inf_{\boldsymbol{\xi}\in
    \widetilde{\Xi}} \textrm{dist}
    (\boldsymbol{\xi}_j,\boldsymbol{\xi}) = \inf_{\boldsymbol{\xi}\in
    \widetilde{\Xi}}||\boldsymbol{\xi}_j-\boldsymbol{\xi}||.
$
\end{definition}

The \emph{minimum transportation budget} is a minimum threshold to ensure the non-emptiness of the trimmings-based ambiguity set, which is finally defined as follows.


 }

%

{\color{black}
 \begin{definition}[Ambiguity set based on probability trimmings]\label{def_ambiguity_set_trimm}
 Assume that $\mathbb{Q}\in \mathcal{P}_1(\Xi)$,  $\alpha>0$, and consider the set $\mathcal{R}_{1-\alpha}(\widehat{\mathbb{Q}}_N)$ of  $(1-\alpha)-$trimmings of the empirical joint distribution $\widehat{\mathbb{Q}}_N:=\frac{1}{N}\sum_{i=1}^N \delta_{(\widehat{\mathbf{z}}_i,\widehat{\boldsymbol{\omega}}_i)}$.
 If $\rho\geqslant\underline{\rho}_{N\alpha}$, where $\underline{\rho}_{N\alpha}$ is the minimum transportation budget introduced in Definition~\ref{MTB},
 then
  the ambiguity set based on trimmings,   $\widehat{\mathcal{U}}_N(\alpha,\rho)$, is defined as the set of all distributions $Q_{\widetilde{\Xi}}$ supported on $\widetilde{\Xi}$ such that $W_1(\mathcal{R}_{1-\alpha}( \widehat{\mathbb{Q}}_N),Q_{\widetilde{\Xi}}) \leq \rho$.
 \end{definition}
 }

 {\color{black} The ambiguity set based on probability trimmings collects all  distributions  with support set $\widetilde{\Xi}$ that result  from a \emph{partial optimal transport} of mass $\alpha$ from $\widehat{\mathbb{Q}}_N$  to $\widetilde{\Xi}$
within a  budget $\rho$. The parameters $\alpha$ and $\rho$ tune the \emph{shape} and the \emph{size} of the ambiguity set, respectively. Indeed,  $\alpha$ can be seen as  the minimum amount of mass of $\widehat{\mathbb{Q}}_N$  which plays a role in the conditional inference and $\rho$ controls the degree of robustness/conservativeness. Hereinafter, we refer to $\alpha$ as the \emph{trimming level} and to $\rho$ as the \emph{transportation budget} or \emph{robustness parameter}.}
%
%
%
Further technical information on the ambiguity set $\widehat{\mathcal{U}}_N(\alpha,\rho)$ can be obtained from~\cite{Estebanmorales_trimmings_extended2020}.


{\color{black}
\begin{remark}
In order to account for contextual information in a DRO framework,  the authors in \cite{Bertsimas2022} propose an ambiguity set $\widehat{\mathcal{U}}$ different to the one based on probability trimmings  (i.e., $\widehat{\mathcal{U}}_N(\alpha,\rho)$)  that we advocate here. More specifically, the ambiguity set they suggest is a Wasserstein ball centered at the discrete distribution supported on the $\boldsymbol{\widehat{\omega}}$-coordinates of the $K$ data points from the sample $\widehat{\Xi}^N_{\omega}$ that are the closest to $\widetilde{\Xi}$. In {\color{black}the supplementary material}, we use an example based on a small three-node system {\color{black}  which  has been taken from \cite{Morales2012}} to illustrate that our trimmings-based ambiguity set generally delivers better dispatch solutions in terms of expected cost and system reliability than the one introduced in \cite{Bertsimas2022} for the DC-OPF problem under uncertainty.
\end{remark}
  }


In the next section, we introduce a tractable and conservative  approximation of the distributionally robust joint chance constraints \eqref{chance_constraintsDRO} using the notion of Conditional Value at Risk (${\mathbf{CVaR}}$). As previously mentioned, the use of the ${\mathbf{CVaR}}$  in the context of chance-constrained distributionally robust OPF is very popular in the  literature (see, for example, the recent publications \cite{Jabr2020} and  \cite{Ordoudis2021}).

\section{A tractable  approximation of the distributionally robust joint chance constraints}\label{sec:CVaRapprox}

The distributionally robust joint chance constraints \eqref{chance_constraintsDRO}  can be written equivalently as the following distributionally robust single chance constraint, where we have already replaced the generic $\widehat{\mathcal{U}}$ with the ambiguity set based on probability trimmings $\widehat{\mathcal{U}}_N(\alpha,\rho)$:
 {\color{black}
 \begin{equation}\label{DRCC-single}
     \inf_{Q_{\widetilde{\Xi}}\in \widehat{\mathcal{U}}_N(\alpha,\rho) }Q_{\widetilde{\Xi}}\left( \max_{k \leq K} \langle \mathbf{a}_{1k},\boldsymbol{\omega} \rangle +a_{2k}\leqslant 0\right)\geqslant 1-\epsilon
 \end{equation}
 Functions $\langle \mathbf{a}_{1k},\boldsymbol{\omega} \rangle +a_{2k}$, $k \leq K$, represent the OPF constraints involved in the joint chance constraint  \eqref{chance_constraintsDRO} } expressed as inequalities lower than or equal to zero. {\color{black} These constraints are all linear with respect to the random vector $\boldsymbol{\omega}$.}


{\color{black} In practice, the system operator is not only concerned about the joint violation of the OPF constraints, but also about the magnitude of this violation. Indeed, the joint chance constraint \eqref{DRCC-single} does not offer guarantees on how positive $ \max_{k \leq K} \langle \mathbf{a}_{1k},\boldsymbol{\omega} \rangle +a_{2k}$ is. This is one of the main reasons to adopt a risk-averse approach to handle the joint chance constraint via the well-known concept of Conditional-Value-at-Risk (${\mathbf{CVaR}}$), which quantifies the conditional expectation   of $ \max_{k \leq K} \langle \mathbf{a}_{1k},\boldsymbol{\omega} \rangle +a_{2k}$ on its right $\epsilon$-tail. More specifically,  the
${\mathbf{CVaR}}$ at level $\epsilon \in (0,1)$ of an univariate random variable $\phi(\boldsymbol{\omega})$
under the probability measure $Q$, $Q-{\mathbf{CVaR}}_{\epsilon}(\phi(\boldsymbol{\omega}))$, is defined as    the value  $\inf_{\tau \in \mathbb{R}}\{\tau +\frac{1}{\epsilon}\mathbb{E}_{Q}[(\phi(\boldsymbol{\omega})-\tau)^{+}] \}$ and when the infimum is attained, $\tau$ represents the Value-at-Risk with confidence level $1-\epsilon$ \cite{Rockafellar00optimizationof}.

}

 In this paper, in lieu of~\eqref{DRCC-single}, we consider the following tractable (convex) approximation:
\begin{equation}\label{DRCC-CVAR}
    \sup_{Q_{\widetilde{\Xi}}\in \widehat{\mathcal{U}}_N(\alpha,\rho) }Q_{\widetilde{\Xi}}-{\mathbf{CVaR}}_{\epsilon}\left( \max_{k \leq K} {\color{black}\langle \mathbf{a}_{1k},\boldsymbol{\omega} \rangle +a_{2k}}\right)\leqslant 0
\end{equation}
which is, in addition, conservative, because \eqref{DRCC-CVAR} implies \eqref{DRCC-single}.

Constraint \eqref{DRCC-CVAR} can be equivalently reformulated as follows \cite{Ordoudis2021}, \cite{Rockafellar00optimizationof}:
\begin{equation}\label{cvar2}
 \inf_{\tau \in \mathbb{R}} \left\{\tau +\frac{1}{\epsilon}\sup_{Q_{\widetilde{\Xi}} \in\widehat{\mathcal{U}}_N(\alpha,\rho) } \mathbb{E}_{Q_{\widetilde{\Xi}}}\left[ \left( \max_{k \leq K} {\color{black}\langle \mathbf{a}_{1k},\boldsymbol{\omega} \rangle +a_{2k}}-\tau\right)^+ \right] \right\} \leqslant 0
\end{equation}
The next proposition states a tractable reformulation of \eqref{cvar2}. For ease of exposition, we first need to recast function $(\max_{k \leq K}\langle \mathbf{a}_{1k},\boldsymbol{\omega} \rangle +a_{2k}-\tau)^+$ as
\begin{equation}
\left(\max_{k \leq K}{\color{black}\langle \mathbf{a}_{1k},\boldsymbol{\omega} \rangle +a_{2k}}-\tau\right)^+:= \max_{k\leqslant K+1}{\color{black} \langle \mathbf{a}_{1k},\boldsymbol{\omega} \rangle +a'_{2k}}\end{equation}
{\color{black}where $a'_{2k}=a_{2k}-\tau$  for  $k\leq K$,  $\mathbf{a}_{1K+1}=\mathbf{0}$ and $a'_{2K+1}=0$.}

\begin{prop}[Reformulation of the ${\mathbf{CVaR}}$-based distributionally robust joint chance constraints]\label{reform_cvar_constraints}
Set $\alpha > 0$.
Then, for any value of $\rho\geqslant\underline{\rho}_{N\alpha}$, the ${\mathbf{CVaR}}$-based distributionally robust joint chance constraints  defined by \eqref{DRCC-CVAR} can be equivalently reformulated  as follows:
\begin{subequations}\label{Prop1_inf}
\begin{align}
\inf_{\tau \in \mathbb{R},\lambda_2 \geqslant 0,\mu_i \geqslant 0, \theta_2 \in \mathbb{R},{\color{black}\boldsymbol{\gamma}_{ik}},\mathbf{v}_{ik}}& \left\{\tau +\frac{1}{\epsilon}\left[ \lambda_2\rho+\theta_2+\frac{1}{N\alpha}\sum_{i=1}^N\mu_i\right] \right\} \leqslant 0    \\
{\text {s.t.}}\,&\  \mu_i+\theta_2 +\lambda_2 \|\mathbf{z}^*-\widehat{\mathbf{z}}_i \|_1 \geqslant {\color{black}a'_{2k}}+S_{\widetilde{\Xi}_{\mathbf{\omega}}}(\mathbf{v}_{ik})\notag \\
&-\langle {\color{black}\boldsymbol{\gamma}_{ik}}, \widehat{\boldsymbol{\omega}}_i\rangle,\forall i \leqslant N, \forall k \leqslant K+1\\
&{\color{black}\boldsymbol{\gamma}_{ik}}-\mathbf{v}_{ik}=-\mathbf{a}_{1k},\forall i \leqslant N, \forall k \leqslant K+1\\
&\|{\color{black}\boldsymbol{\gamma}_{ik}} \|_{\infty}\leqslant \lambda_2 ,\forall i \leqslant N, \forall k \leqslant K+1
\end{align}
\end{subequations}
where $S_{\widetilde{\Xi}_{\mathbf{\omega}}}(\cdot)$ stands for the support function of $\widetilde{\Xi}_{\mathbf{\omega}}$ and $\langle \cdot, \cdot\rangle$ represents the dot product (see \ref{appendix:othersymbols}).

\end{prop}


Once we have reformulated the ${\mathbf{CVaR}}$-based distributionally robust joint chance constraint~\eqref{DRCC-CVAR}, we only need to reformulate the DRO problem defined by the inner supremum in \eqref{worst_case_expected}. Since this requires a careful and independent analysis, we consider it in the following section.
\section{An exact tractable reformulation of the worst-case expected  cost}\label{sec:worstEC}
In what follows, we provide an exact and tractable reformulation of the objective function \eqref{worst_case_expected} as a continuous linear program.

The term {\color{black}
\begin{equation}\label{reform_max_piece1}
 C(\tilde{\mathbf{g}}(\boldsymbol{\omega}))=   \sum_{j\in \mathcal{G}} \max_{s=1,\ldots, S_j } \Bigg\{ m_{js}\left[g_j-\beta_j\Omega\right] +n_{js}\Bigg\}
\end{equation}
}
is {\color{black}the} sum of {\color{black}the} maximum of \emph{univariate} linear functions in terms of $\Omega$, which is, moreover, convex in $\Omega$. This observation is key to reformulating \eqref{worst_case_expected} in a tractable way. In fact, the ambiguity set $\widehat{\mathcal{U}}_N(\alpha,\rho)$ for the worst-case probability distribution in the inner supremum of \eqref{worst_case_expected} can  be equivalently replaced with the following one, which is also expressed in terms of $\Omega$ only:
\begin{equation}\label{ambiguity_set_Omega}
    \widehat{\mathcal{U}}_N^{\Omega}(\alpha, \rho):=\{P_{\widetilde{\Xi}_{\Omega}} : W_1(\mathcal{R}_{1-\alpha}( \widehat{\mathbb{P}}_N),P_{\widetilde{\Xi}_{\Omega}}) \leq \rho,\; P_{\widetilde{\Xi}_{\Omega}}(\widetilde{\Xi}_{\Omega})=1  \}
\end{equation}
where $\widehat{\mathbb{P}}_N:=\frac{1}{N}\sum_{i=1}^N \delta_{(\widehat{\mathbf{z}}_i, \widehat{\Omega}_i)}$ is the empirical distribution supported on the samples $(\widehat{\mathbf{z}}, \widehat{\Omega}_i), i=1,\ldots,N$; and $\widetilde{\Xi}_{\Omega}$ stands for the event  $$(\mathbf{z}=\mathbf{f}; \Omega \in [\underline{\Omega},\overline{\Omega}]), \; \text{with}\; [\underline{\Omega},\overline{\Omega}]=\left[-\sum_{m\in \mathcal{W}}f_m,\sum_{m\in \mathcal{W}}(\overline{C}_m-f_m)\right]$$
The interval $[\underline{\Omega},\overline{\Omega}]$ is the conditional support for the random variable $\Omega$ (that is, the support set for the system-wise aggregate wind power forecast error{\color{black},} given the predicted power outputs of the wind farms). Essentially, what we have done above is to map the original probability space for the random vector $(\mathbf{z},\boldsymbol{\omega})$ onto a new probability space for the random vector $(\mathbf{z},\Omega)$ by the {\color{black}linear map $\boldsymbol{\omega}\mapsto\sum_{m\in \mathcal{W}}\omega_m$, $\Omega=\sum_{m\in \mathcal{W}}\omega_m$}, which leaves the objective cost function unaltered. In doing so, the inner supremum in \eqref{worst_case_expected} can be fully recast in terms of $\Omega$ only as follows:
\begin{equation}\label{DRO_Omega_objvalue}
    \sup_{P_{\widetilde{\Xi}_{\Omega} }\in \widehat{\mathcal{U}}^{\Omega}_N(\alpha,\rho)} \mathbb{E}_{P_{\widetilde{\Xi}_{\Omega} }} \left[   {\color{black}  \sum_{j\in \mathcal{G}} \max_{s=1,\ldots, S_j } \Bigg\{ m_{js}\left[g_j-\beta_j\Omega\right] +n_{js}\Bigg\}+\langle \mathbf{c}^D, \mathbf{r}^D \rangle +\langle \mathbf{c}^U, \mathbf{r}^U \rangle } \right]
\end{equation}

The proposition below presents a tractable reformulation of \eqref{DRO_Omega_objvalue} as a continuous linear program.

\begin{prop}[LP reduction of the worst-case expected cost]\label{lpreduction_objective}
Set $\alpha > 0$ and assume that  $\|(\mathbf{z},\Omega)\|:=\| \mathbf{z} \|+|\Omega|$ for some norm $\| \cdot \|$ in $\mathbb{R}^{d_{\mathbf{z}}}$.  Then, for any value of $\rho\geqslant\underline{\rho}_{N\alpha}$, the DRO problem defined by \eqref{DRO_Omega_objvalue} can be reformulated as the following continuous linear program:
\begin{subequations}
  \begin{align}
  &  \inf_{\lambda \geqslant 0;\theta \in \mathbb{R}, \overline{\mu}_i, {\color{black}t_i}, {\color{black}\underline{t}_{ij},\overline{t}_{ij},\widehat{t}_{ij}  \forall i\leqslant N,\; \forall j \in \mathcal{G}}   } \quad   \lambda \rho+\theta +\dfrac{1}{N\alpha}\sum_{i=1}^N \overline{\mu}_i +{\color{black}\langle \mathbf{c}^D, \mathbf{r}^D \rangle +\langle \mathbf{c}^U, \mathbf{r}^U \rangle}  \\
& \hspace{2cm}  \text{s.t.} \ \overline{\mu}_i+\theta+\lambda\|\mathbf{z}^*-\widehat{\mathbf{z}}_{i}\| \geqslant \;
t_i, \forall i \leqslant N\\
& \hspace{2cm} \phantom{s.t. } t_i \geqslant  \sum_{j\in \mathcal{G}} \underline{t}_{ij} -\lambda (\underline{\Omega}-\widehat{\Omega}_i), \; \forall i \in \underline{I}\\
& \hspace{2cm}\phantom{s.t. }  t_i \geqslant  \sum_{j\in \mathcal{G}} \overline{t}_{ij} -\lambda (\overline{\Omega}-\widehat{\Omega}_i), \; \forall i \in \underline{I}\\
& \hspace{2cm}\phantom{s.t. }  t_i \geqslant  \sum_{j\in \mathcal{G}}\underline{t}_{ij} +\lambda (\underline{\Omega}-\widehat{\Omega}_i), \; \forall i \in \overline{I}\\
& \hspace{2cm}\phantom{s.t. }  t_i \geqslant  \sum_{j\in \mathcal{G}} \overline{t}_{ij} +\lambda (\overline{\Omega}-\widehat{\Omega}_i), \; \forall i \in \overline{I}\\
& \hspace{2cm}\phantom{s.t. }  t_i \geqslant  \sum_{j\in \mathcal{G}}\overline{t}_{ij} -\lambda (\overline{\Omega}-\widehat{\Omega}_i), \; \forall i \in I\\
& \hspace{2cm} \phantom{s.t. } t_i \geqslant  \sum_{j\in \mathcal{G}} \underline{t}_{ij} +\lambda (\underline{\Omega}-\widehat{\Omega}_i), \; \forall i \in I\\
& \hspace{2cm} \phantom{s.t. } t_i \geqslant  \sum_{j\in \mathcal{G}} \widehat{t}_{ij} , \; \forall i \in I\\
& \hspace{2cm}  \phantom{s.t. }  \overline{\mu}_i \geqslant 0,\enskip \forall i \leqslant N\\
& \hspace{2cm}\phantom{s.t. } \underline{t}_{ij} \geqslant  m_{js}\left[g_j-\beta_j\underline{\Omega}\right] +n_{js}, \enskip \forall i \leqslant N, \enskip \forall j \in \mathcal{G}, \enskip \forall s \leqslant S_j\\
& \hspace{2cm} \phantom{s.t. }\overline{t}_{ij} \geqslant  m_{js}\left[g_j-\beta_j\overline{\Omega}\right] +n_{js}, \enskip \forall i \leqslant N, \enskip \forall j \in \mathcal{G}, \enskip \forall s \leqslant S_j\\
& \hspace{2cm}\phantom{s.t. } \widehat{t}_{ij} \geqslant  m_{js}\left[g_j-\beta_j\widehat{\Omega}_i\right] +n_{js}, \enskip \forall i \leqslant N, \enskip \forall j \in \mathcal{G}, \enskip \forall s \leqslant S_j
\end{align}
\end{subequations}
where $
   \underline{I}:=\{ i \in \{1, \ldots N\} : \widehat{\Omega}_i < \underline{\Omega} \},
   I:=\{ i \in \{1, \ldots N\} : \widehat{\Omega}_i \in [\underline{\Omega},\overline{\Omega}] \}$, and $\overline{I}:=\{ i \in \{1, \ldots N\} : \widehat{\Omega}_i >\overline{\Omega} \}
$.

\end{prop}

\section{Numerical results}\label{sec:numerics}
In this section, we present and discuss results from a series of numerical experiments that have been run on a modified version of the IEEE 118-bus system considered in \cite{Jabr2020}. All the data and codes needed to reproduce those experiments are available for download in the GITHUB repository \cite{DRO_DCOPF_CONTEXTUAL_github2021}. The experiments have been carried out on a Linux-based server {\color{black} using up to 13200 CPUs running in paralell, each clocking at 2.6 GHz with 200 GB of RAM. We have employed CPLEX 20.1.0 under DOcplex Python Modeling API to solve the associated continuous linear programs with the barrier algorithm using up to 22 threads. In addition, we have set the CPLEX parameter   {\texttt{preprocessing\_dual}} to 1}.

We solve the CC-DRO OPF problem~\eqref{worst_case_expected}--\eqref{chance_constraintsDRO} using the ${\mathbf{CVaR}}$-based approximation stated {\color{black}in Section~\ref{sec:CVaRapprox}}, but with different ambiguity sets, namely: (i) The ambiguity set based on probability trimmings, introduced in this paper, which we refer to as DROTRIMM; and (ii) {\color{black} a Wasserstein ball centered at the empirical distribution supported on the $\widehat{\boldsymbol{\omega}}$-coordinates of the $N$ samples in $\widehat{\Xi}^N_{\omega}$, {\color{black} i.e., $\{(\widehat{\mathbf{z}}_i, \widehat{\boldsymbol{\omega}}_i)\}_{i=1}^N$}. This leads to the distributionally robust chance-constrained OPF model proposed in~\cite{Ordoudis2021}, which we call DROW. Importantly, this is a DRO model that fully ignores the contextual information, since the center of the Wasserstein ball it uses is made up of \emph{all} past samples of wind power forecast errors (regardless of the current wind power point predictions)}.
Roughly speaking, DROTRIMM also works with all the past $N$ samples of wind power forecast errors, but only  those that lead to the worst-case conditional distribution of the prediction errors are moved onto the conditional support. However, this movement must entail a transportation cost  smaller than or equal to a given budget ${\color{black}\rho}$ and the computation of that cost is directly contingent on the current context (that is, the current wind power point forecasts).

{\color{black} In addition, we benchmark the previous two distributionally robust methods with an alternative approach that is commonly used in the  literature for solving optimization problems with probabilistic constraints, known as the \emph{scenario approach}, but adapted to account for contextual information. We have taken the required adaptation from \cite[Chapter 4]{koduri2021essays}, which, in our setting, involves solving a DC-OPF problem in which the uncertain constraints are enforced for the wind power forecast errors associated with the $K$ samples nearest to the context. We refer to this adaptation of the popular scenario approach as SCENA.}

The training data consist of a set of pairs $\{(\widehat{\mathbf{z}}_i, \widehat{\boldsymbol{\omega}}_i)\}_{i=1}^N$, from which we can directly obtain the collection of pairs $\{   (\widehat{\mathbf{z}}_i, \widehat{\Omega}_i)\}_{i=1}^N$, where $ \widehat{\Omega}_i:=\sum_{{\color{black}m\in \mathcal{W}}} \widehat{\omega}_{i,m}$. For ease of computation and to simplify the analysis below,  we have considered the same radius or transportation budget for the two ambiguity sets in both the objective and the chance constraints of the DRO OPF problem~{\color{black} \eqref{worst_case_expected}--\eqref{chance_constraintsDRO}}.

\subsection{Evaluation of the out-of-sample performance via re-optimization}

Given a context (in the form of point forecasts of the power outputs of the wind farms), a training dataset, and a robustness parameter $\rho$, each method $met$ (either {\color{black} DROW} or DROTRIMM in our case) provides a forward generation dispatch and reserve capacity provision $\mathbf{y}^{met}:=(\mathbf{g},\mathbf{r}^D,\mathbf{r}^U)$. To evaluate the actual or out-of-sample performance of that $\mathbf{y}^{met}$, we draw a sample of wind power forecast errors $\widehat{\boldsymbol{\omega}}_j$ from a test dataset, and the vector of recourse variables
$\mathbf{r}$ (that is, the real-time power adjustments) is computed by solving the deterministic optimal power flow available in \ref{appendix_real_time_formulation}. In this deterministic OPF problem, wind spillage (with a cost equal to 0) and involuntary load curtailment (with a cost equal to \$500/MWh) are considered as feasible recourse actions, aside from the deployment of reserves by generators.
In this way, the \emph{out-of-sample performance} of  a method $met$, which produces the  forward dispatch $\mathbf{y}^{met}$, $J(\mathbf{y}^{met})$, is computed by the empirical out-of-sample  cost {\color{black} averaged} over the test set formed by a certain number of samples of $\mathbb{Q}_{\boldsymbol{\omega}/\mathbf{z}=\mathbf{f}}$. In addition, in order to  measure the \emph{reliability} of a solution (that is, if $\mathbf{y}^{met}$ is feasible or not in real time), the \emph{violation probability} is estimated over the  test set. {\color{black} In this estimation, we count as a violation every time a recourse action involving load curtailment or wind spillage is to be taken in real time to restore the power balance. This is equivalent to counting (over the test set) the number of times a constraint is violated in the original affine-policy-based OPF model.}

\subsection{A 118-bus case study}
As previously mentioned, we consider a modified version of the IEEE 118-bus system used in \cite{Jabr2020}. The system includes 54 conventional generators and {\color{black} eight wind power plants that we have added and placed at buses  2, 16, 33, 37, 55, 67, 83, and  116}. In addition, the piecewise linear cost functions of all generators are comprised of three pieces or blocks. All  the data pertaining to the  network, generators,  and  transmission lines are available at the GITHUB repository \cite{DRO_DCOPF_CONTEXTUAL_github2021}.

We analyze two scenarios, which differ in the level of wind power penetration in the system. Below we explain how we have generated samples for the joint distribution of the wind power forecast  and its error at each wind power plant. The so generated samples are also available online at the GITHUB repository  \cite{DRO_DCOPF_CONTEXTUAL_github2021}:

\begin{enumerate}
    \item Let $\widetilde{f}_m$ be the per-unit point forecast of the power output  at wind plant ${\color{black}m\in \mathcal{W}}$. {\color{black} A sample of $\widetilde{f}_m$, for all $m\in \mathcal{W}$, is randomly drawn with replacement from a collection of  16\,694 p.u. wind power data} recorded in several zones and made available by the Global Energy Forecasting Competition 2014 \cite{global_forecast_competition2014}. {\color{black} We have selected zones 1, 2, 3, 4, 5, 6, 9, and 10 of the aforementioned data set and assigned them to the eight wind power plants located at buses 16, 116, 83, 2, 55, 67, 33, and 37,  respectively.}

    \item For each wind farm ${\color{black} m \in \mathcal{W}}$, we have assumed that the  (nominal, normalized) random variable $W_m$, which represents the nominal actual power generated at wind plant $m$, follows a Beta distribution with mean $\widetilde{f}_m$ and standard deviation $\sigma$. This standard deviation depends on both physical parameters and the quality of the forecasting model, following  the model proposed in \cite{Fabbri2005}. For simplicity, in all numerical experiments, given $\widetilde{f}_m$,  we determine $\sigma$ as the value of the following function  $\sigma(\widetilde{f}_m):=0.2\widetilde{f}_m+0.02$,
    empirically obtained in \cite{Fabbri2005} for the case of a lead time of six hours. Therefore, the actual wind power production, and hence, the forecast error are conditional on the  forecast power output issued. More specifically, the forecast error is given as the difference of a realization $\widehat{W}_{m}$ of the r.v. $W_m \sim   \text{Beta}(A,B)$ and the point forecast $\widetilde{f}_m$, where $A,B>0$ are the solution (if it exists) of the following system of non-linear equations: \begin{subequations}
    \begin{align}
       \widetilde{f}_m &=\frac{A}{A+B} \\
      \sigma^2(\widetilde{f}_m)  &=\frac{AB}{(A+B)^2(A+B+1)}
    \end{align}
    \end{subequations}
    To ensure that this non-linear system of equations has a solution, the samples $\widetilde{f}_m$ from the dataset that are less than or equal to 0.05 p.u. are set to 0.05, and the ones greater than or equal to 0.95 p.u. are set to 0.95. Thus, for each wind power plant, the per-unit point forecast lies in the interval $[0.05,0.95]$.
     \item To work with MW, we multiply the per-unit realized power output and the point prediction $\widetilde{f}_m$ by the wind plant capacity $\overline{C}_m$, thus getting a pair of predicted power output and its error $(\overline{C}_m\widetilde{f}_m, \overline{C}_m(\widehat{W}_{m}-\widetilde{f}_m))$ for wind farm $m$.
     \item {\color{black} Steps 1, 2 and 3 are repeated $N$ times to get the desired sample size.}
\end{enumerate}

{\color{black} Each independent run in our simulations involves repeating the above process.

Finally, the test set used to compute the out-of-sample performance of a data-driven solution via re-optimization (i.e., the actual probability of violating the uncertain OPF constraints and the actual expected operational cost) is constructed by drawing 1000 samples from
the wind-power-data generating model based on the beta distribution presented above, with mean equal to the point prediction acting as the selected context. Therefore,} this test set constitutes a discrete approximation of the forecast error distribution conditional on a given context, which will be specified later. {\color{black} Importantly, the shape and size of the ambiguity set that DROTRIMM uses is to be changed with the sample size $N$ (which is indicative of the amount of information on the joint distribution of $(\mathbf{z},\boldsymbol{\omega})$ we have). Consequently, the trimming level $\alpha$ defining this set is to be dependent on $N$. Accordingly, we have set $\alpha_N:=K_N/N$, where $K_N$ is the number of nearest neighbors used by SCENA. We have specifically taken $K_N:=\lfloor N^{0.9}\rfloor$ so that the resulting $\alpha_N$ is consistent with the convergence results included in \cite{Estebanmorales_trimmings_extended2020}. Observe that both $\alpha$ and $K$ have been augmented with the subscript $N$ to make their dependence on the sample size explicit.}

\subsubsection{Medium wind penetration case}
In this case, all  {\color{black}eight wind farms in the system have a capacity of 200 MW and the context is given by  $\mathbf{z}^*=180\cdot\mathbf{1}$ MW, that is, the point forecast is 180 MW} for all the wind power plants. Hence, the level of wind power penetration in the system (i.e., the ratio of the predicted system-wise wind power production to the total system demand)  is approximately $63 \%$.

\begin{figure}
\centering
\subfloat[Total downward reserves]{%
  \includegraphics[width=0.5\textwidth]{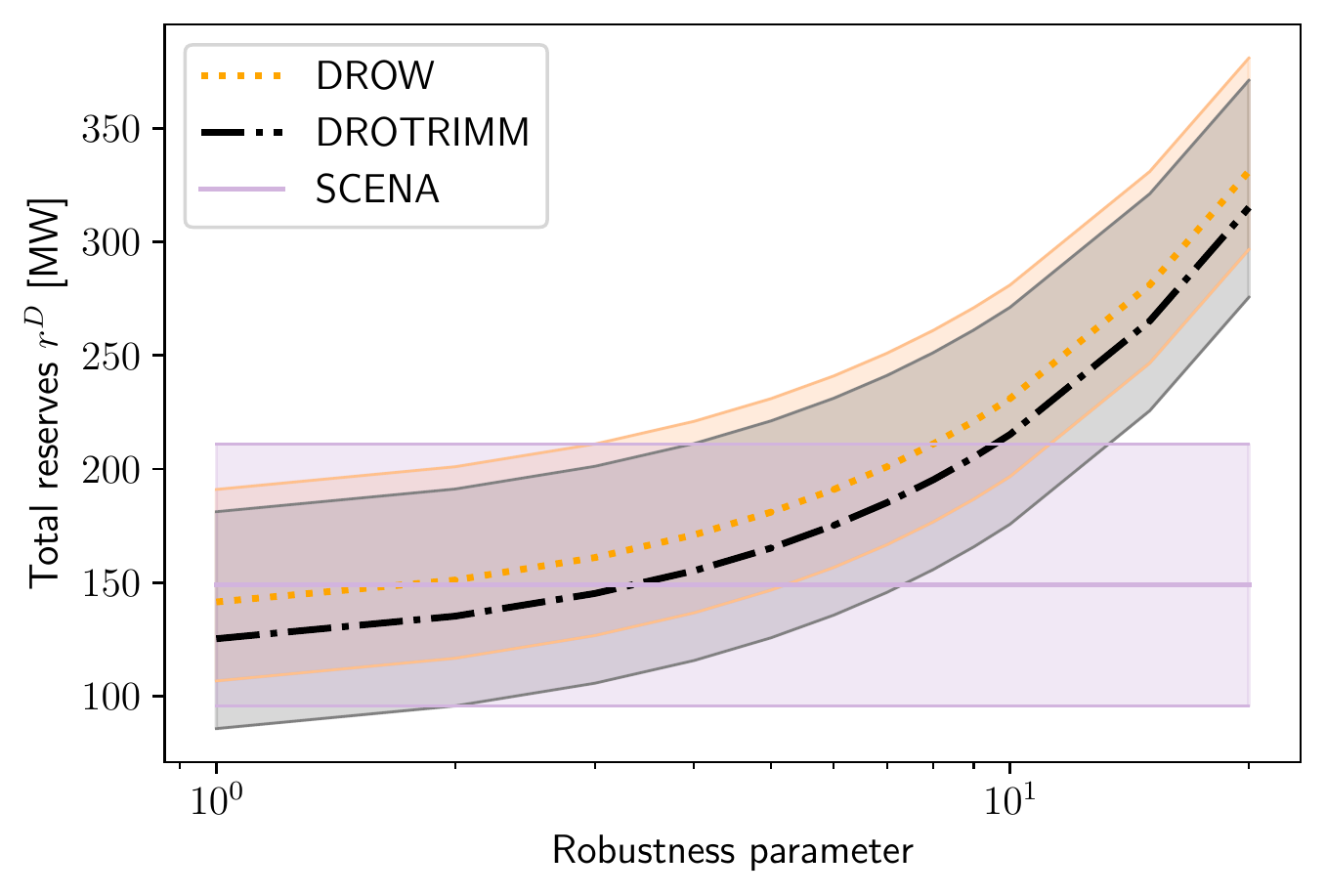}%
  \label{totalreserves_dcopf_sens_100_case_1pos_D}
}%
\subfloat[Total upward reserves]{%
  \includegraphics[width=0.5\textwidth]{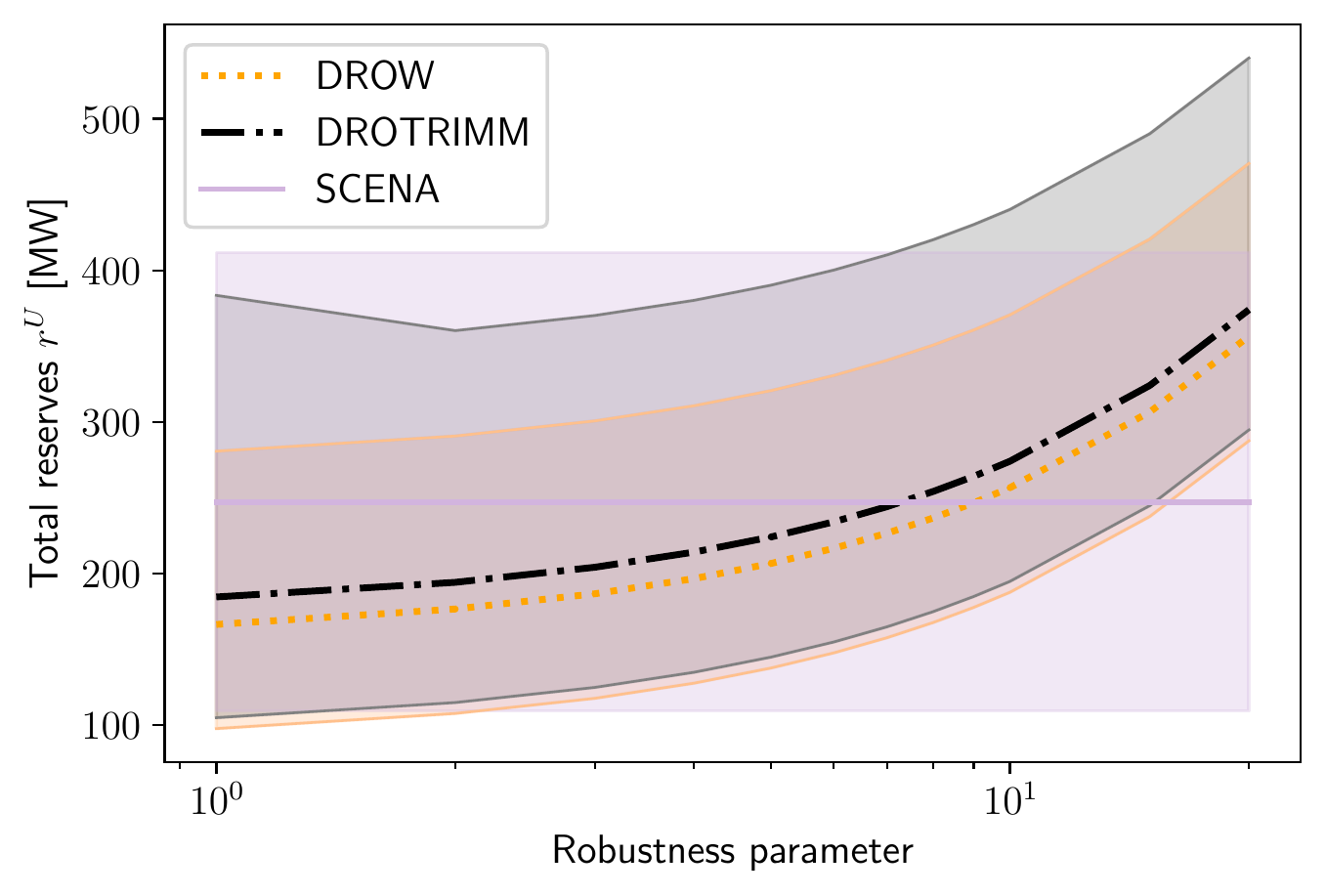}%
  \label{totalreserves_dcopf_sens_100_case_1pos_U}
}

\subfloat[Violation probability (out of sample)]{%
  \includegraphics[width=0.5\textwidth]{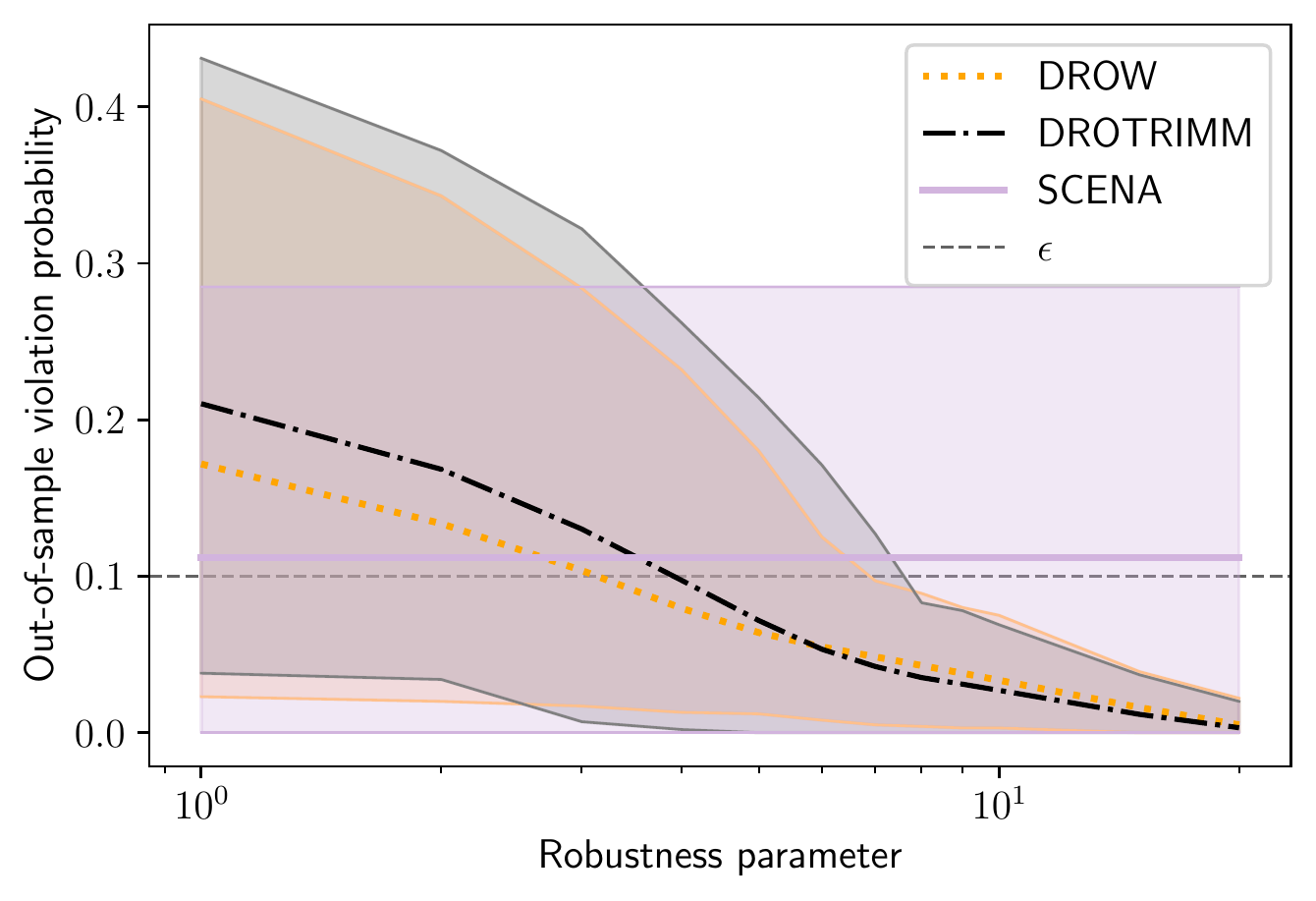}%
  \label{violation_dcopf_sens_100_case1}
}%
\subfloat[Expected system operating cost (out of sample)]{%
  \includegraphics[width=0.5\textwidth]{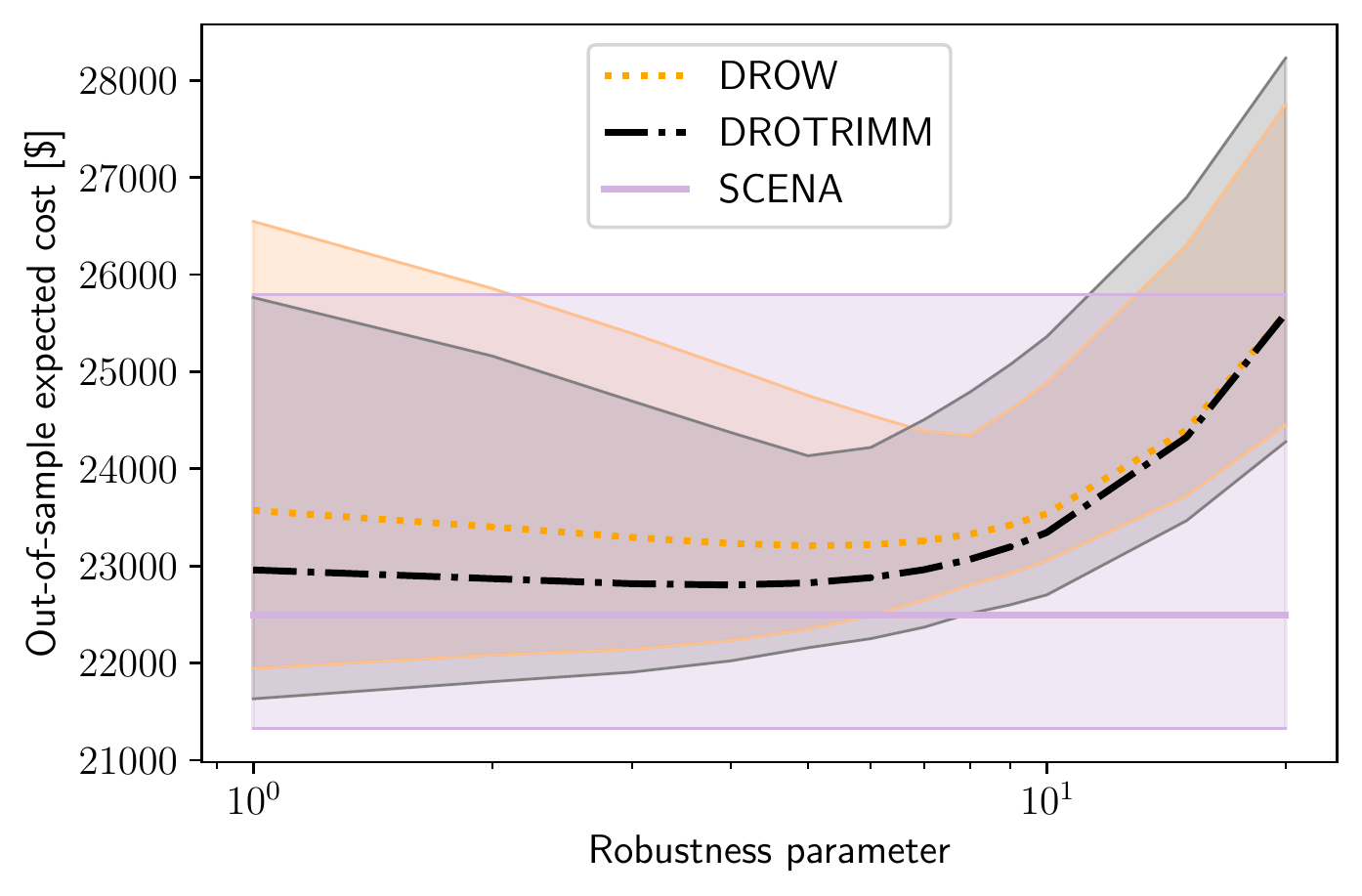}%
  \label{actualexpectedcost_dcopf_sens_100_case1}
}


\caption{Medium level of wind penetration, $N=100 $ and $\epsilon=0.1$:  Total downward and upward reserve capacity and performance metrics}\label{performance_dcopf_118node_100_case1}

\end{figure}

\begin{figure}
\centering
\subfloat[Total downward reserves]{%
  \includegraphics[width=0.5\textwidth]{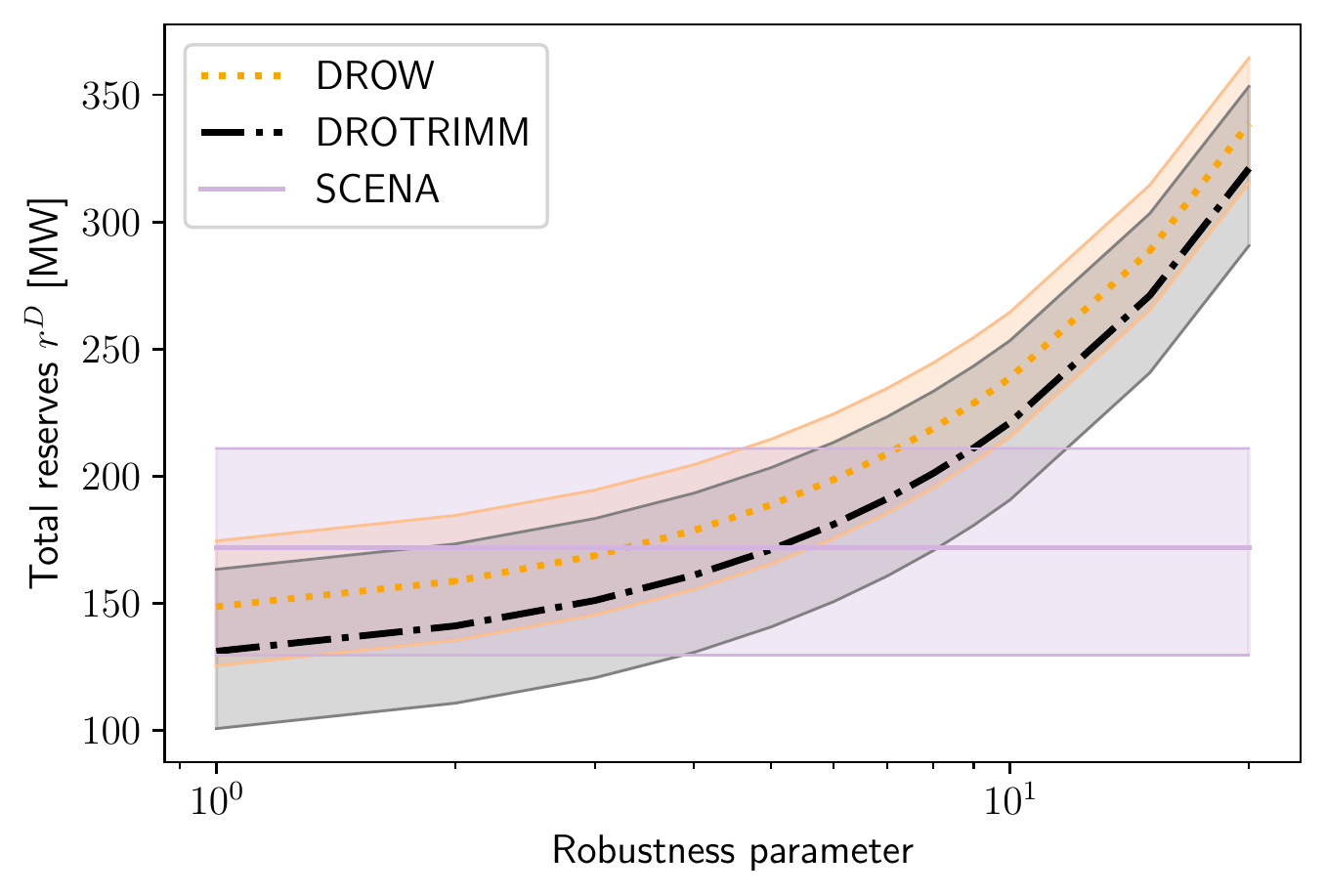}%
  \label{totalreserves_dcopf_sens_300_case_1pos_D}
}%
\subfloat[Total upward reserves]{%
  \includegraphics[width=0.5\textwidth]{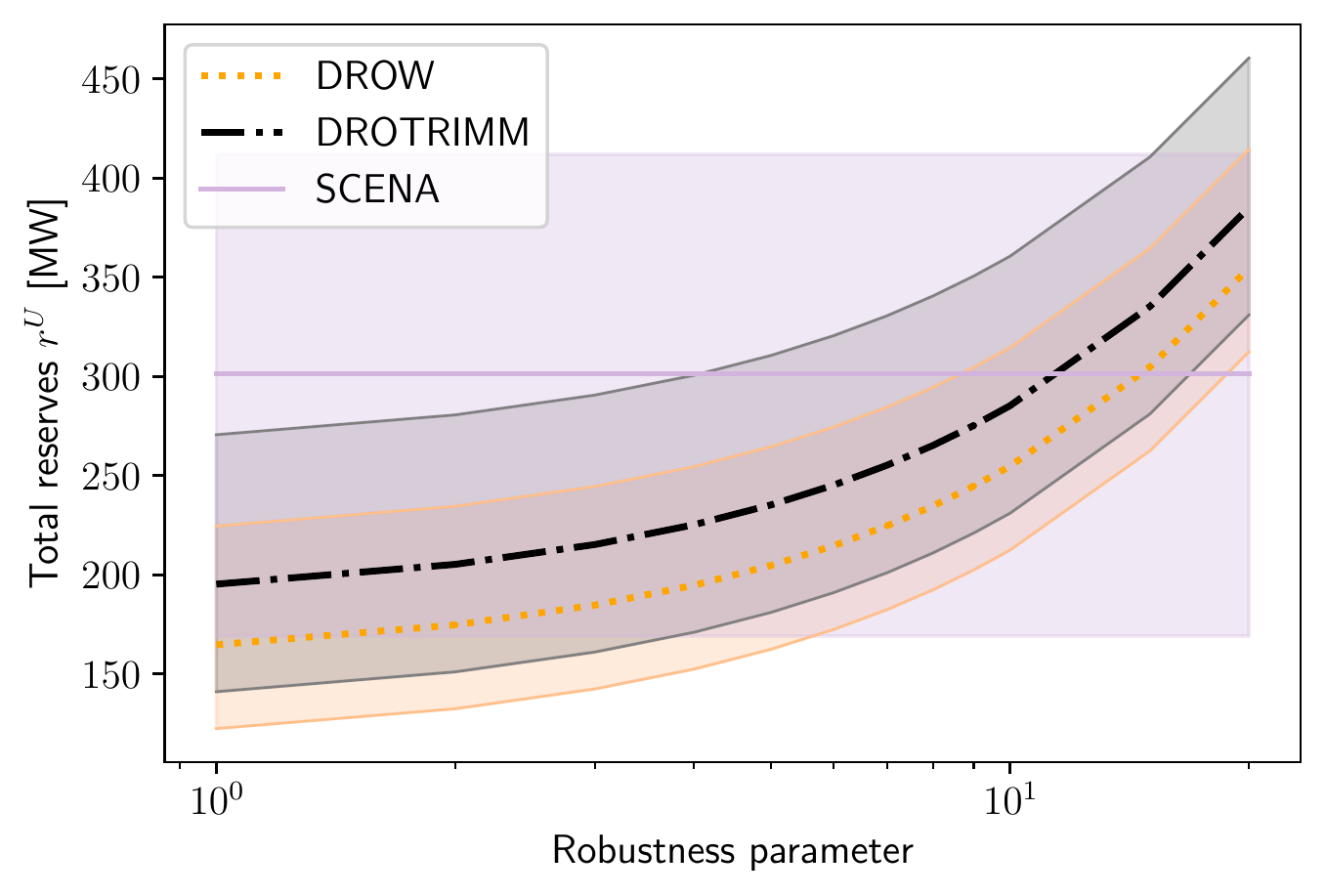}%
  \label{totalreserves_dcopf_sens_300_case_1pos_U}
}

\subfloat[Violation probability (out of sample)]{%
  \includegraphics[width=0.5\textwidth]{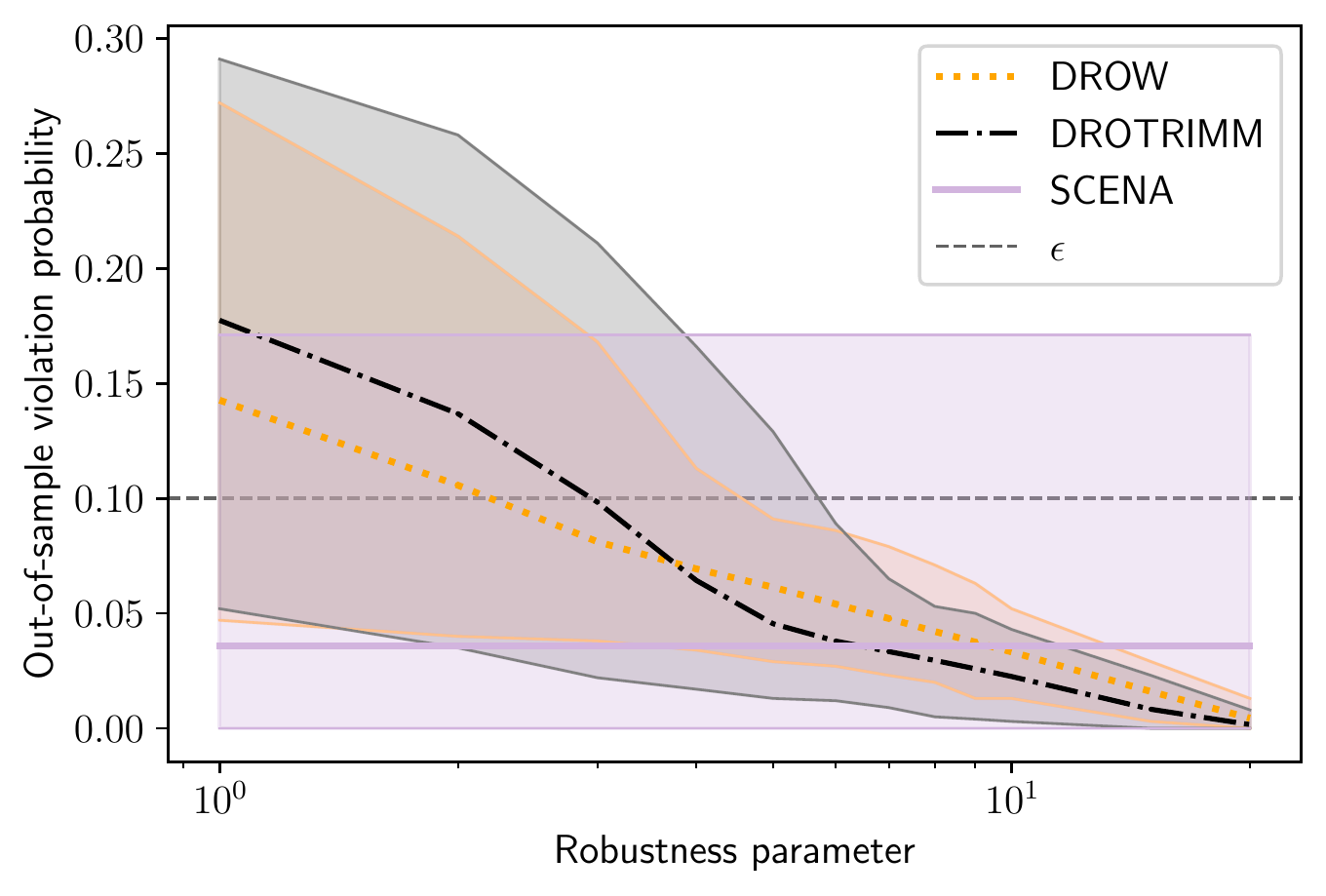}%
  \label{violation_dcopf_sens_300_case1}
}%
\subfloat[Expected system operating cost (out of sample)]{%
  \includegraphics[width=0.5\textwidth]{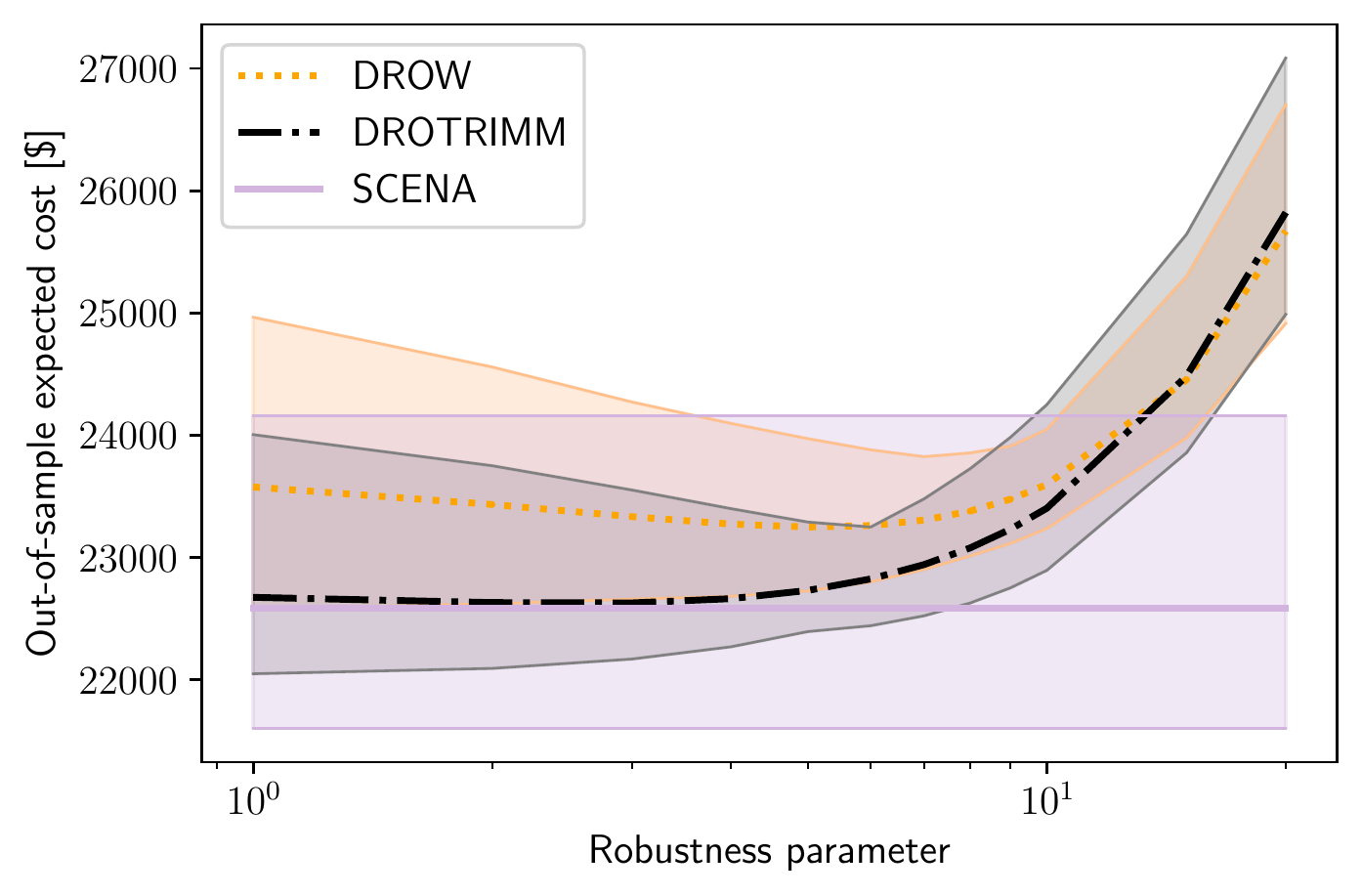}%
  \label{actualexpectedcost_dcopf_sens_300_case1}
}

\caption{Medium level of wind penetration, $N=300 $ and $\epsilon=0.1$:  Total downward and upward reserve capacity and performance metrics}\label{performance_dcopf_118node_300_case1}

\end{figure}

Figures  \ref{performance_dcopf_118node_100_case1} and \ref{performance_dcopf_118node_300_case1}  illustrate the  box
plots  corresponding  to  the  total downward and upward reserve capacity that is scheduled,  the  violation probability and the expected cost delivered out of sample by
{\color{black}SCENA, DROTRIMM and DROW as a function of their corresponding robustness parameter for sample sizes $N=100$ and $N =300$, respectively. Naturally, the results provided by SCENA do not change along the $x$-axis in the plots, because this method is not based on \emph{distributional robustness}}. The box plots have been obtained from 200 independent runs for each sample size. We have set $\epsilon = 0.1$.  The robustness parameter of DROW is the radius of the  Wasserstein ball, while the robustness parameter for DROTRIMM is the budget excess over the minimum transportation budget (see Definition \ref{def_ambiguity_set_trimm}) \cite{Estebanmorales_trimmings_extended2020}.

The  color-shaded
areas have been obtained by joining the {\color{black}minimum and maximum edge cases} of the box plots, while the associated bold colored lines link their means. These figures allow us to check which of the methods provides the most cost-efficient dispatch solutions on average without exceeding the threshold $\epsilon$. As expected, the reliability of the OPF solution given by DROW and DROTRIMM increases as the value of their robustness parameter is augmented, because more reserve capacity is procured. In turn, as more reserve capacity is scheduled, the magnitude and frequency of expensive load shedding events tend to diminish, which explains why the expected system operating cost may also decrease with the robustness parameter. This justifies the use of Distributionally Robust Optimization to tackle the chance-constrained OPF problem. However, when said parameter reaches a large enough value, the expected cost starts to grow quickly, because the cost of procuring additional reserve capacity no longer compensates for the cost savings entailed by the reduction in the amount of curtailed load.

{\color{black}While SCENA provides OPF solutions that are competitive in terms of expected cost, these solutions do not comply with the specified reliability threshold in many of the runs. In addition, the performance of the OPF solutions obtained from SCENA exhibit a high variability, which is clearly due to the fact that this method is a non-robust approach and as such, is highly negatively affected by the uncertainty associated with the conditional inference it must perform.

On the other hand, when comparing DROW and DROTRIMM, whereas the former needs a lower value of the robustness parameter to attain the desired level of solution reliability, DROTRIMM gets to identify OPF solutions that are also reliable, while systematically cheaper on average}. This phenomenon becomes even more evident when we increase the sample size $N$ {\color{black}from 100 to 300}.
Indeed, a richer joint data sample contains more information on the statistical dependence of the wind power forecast error on the associated point prediction, which our approach manages to take advantage of. To give some numbers, if we just consider the range of values for the robustness parameters for which the violation probability is kept below the tolerance~$\epsilon$, {\color{black} the average expected cost savings of DROTRIMM with respect to DROW go from 0.82\%, when $N = 100$ to 1.82\%, when $N = 300$}. {\color{black} From a technical point of view, DROW tends to produce OPF solutions with a higher cost because it underestimates the amount of upward reserve capacity that should be procured, clearly because this method is oblivious to the context and therefore, plans for the \emph{marginal} distribution of the wind power forecast errors and not for the conditional one.}

{\color{black} To elaborate further on the differences among the three methods, Table~\ref{tab:medium} includes the maximum, average, and minimum out-of-sample expected cost\footnote{These statistics are computed over the 200 hundred independent runs.} under the value of the robustness parameter that is \emph{optimal} for methods DROW and DROTRIMM, i.e., which leads to reliable OPF solutions with the minimum average expected cost for each of these two approaches. The standard deviation of this cost is also provided in the last row of Table~\ref{tab:medium}. When $N = 100$, the (exacerbated) robustness of DROW produces OPF solutions with low average cost and variance, although DROTRIMM manages to find OPF solutions that are more economical in expectation. When $N=300$, DROTRIMM clearly beats DROW on all metrics, because the excessive robustness of DROW (which is the result of ignoring the context) no longer pays off. Again, SCENA provides the cheapest OPF solutions on average, but these are useless because they do not satisfy the reliability requirement.}

\begin{table}[t]
\caption{Medium level of wind penetration, summary data for total expected cost [$\$$] under the optimal value of the robustness parameter for methods DROW and DROTRIMM. }
\label{tab:medium}
\centering
\setlength\tabcolsep{3pt} 
{\color{black}\begin{tabular}{ccccccccc}
\hline
\multirow{2}{*}{} &  &\multicolumn{3}{c}{$N=100$}& &\multicolumn{3}{c}{$N=300$}\\
\cline{3-5} \cline{7-9}
 &    &DROTRIMM &DROW &SCENA &&DROTRIMM &DROW &SCENA \\
\hline
    &  max & 24790  &  24388 &  25795  & &23250 &  23972&  24160 \\
           &avg & 23068      &   23258 &  22493 & & 22826 &23250  & 22588 \\
             &min  & 22511      &  22649&   21325 & & 22443 & 22728  &21602  \\
         &std &   300    & 295  &  742& & 159&  246 & 496  \\
\hline
\end{tabular}}
\end{table}


\subsubsection{High wind penetration case}
In this alternative setting, all the wind farms have a capacity of {\color{black} 250 MW and the context is given by  $\mathbf{z}^*=225\cdot\mathbf{1}$} MW. Hence, the level of wind power penetration in the system is approximately $80 \%$.

\begin{figure}
\centering
\subfloat[Total downward reserves]{%
  \includegraphics[width=0.5\textwidth]{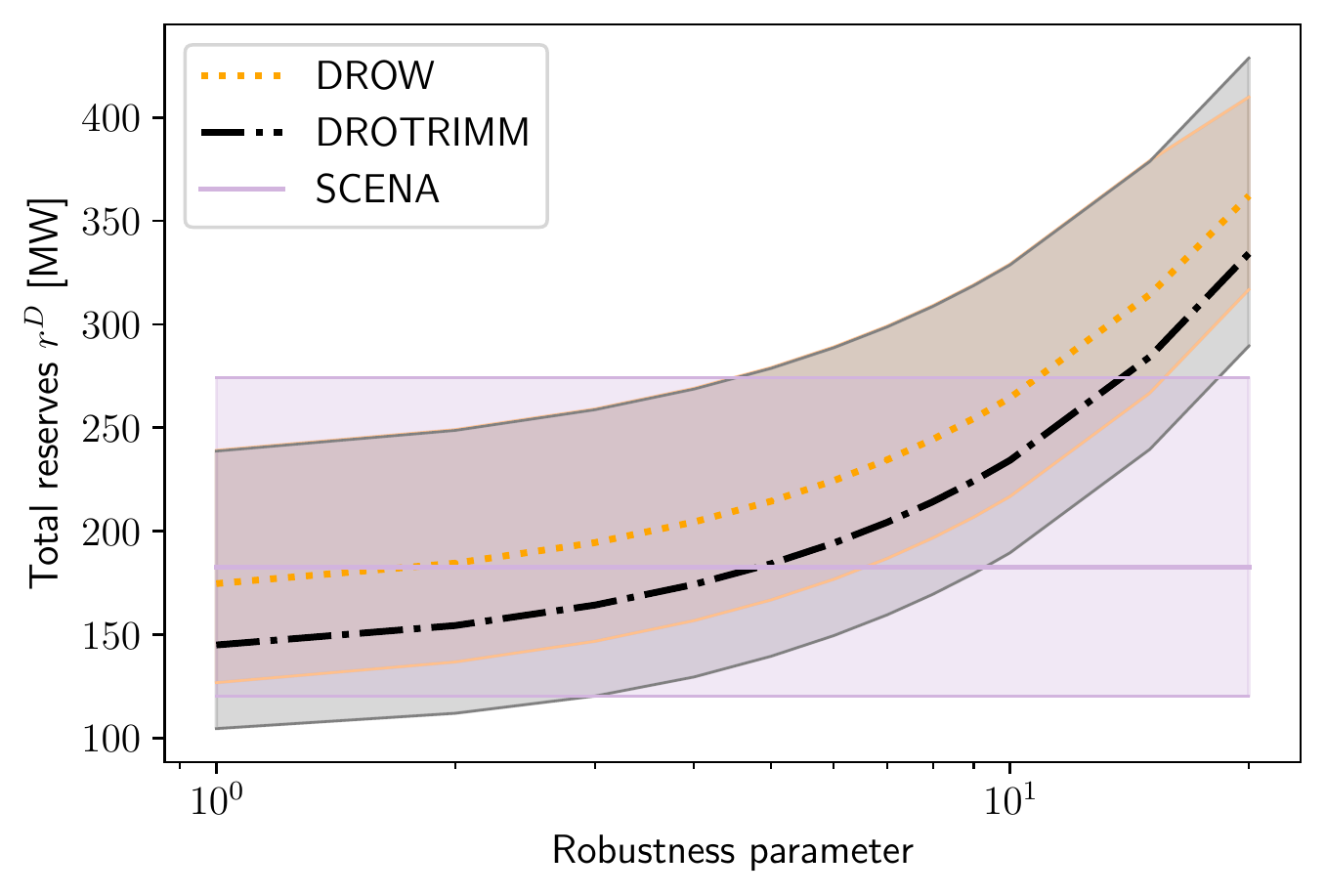}%
  \label{totalreserves_dcopf_sens_100_case_2pos_D}
}%
\subfloat[Total upward reserves]{%
  \includegraphics[width=0.5\textwidth]{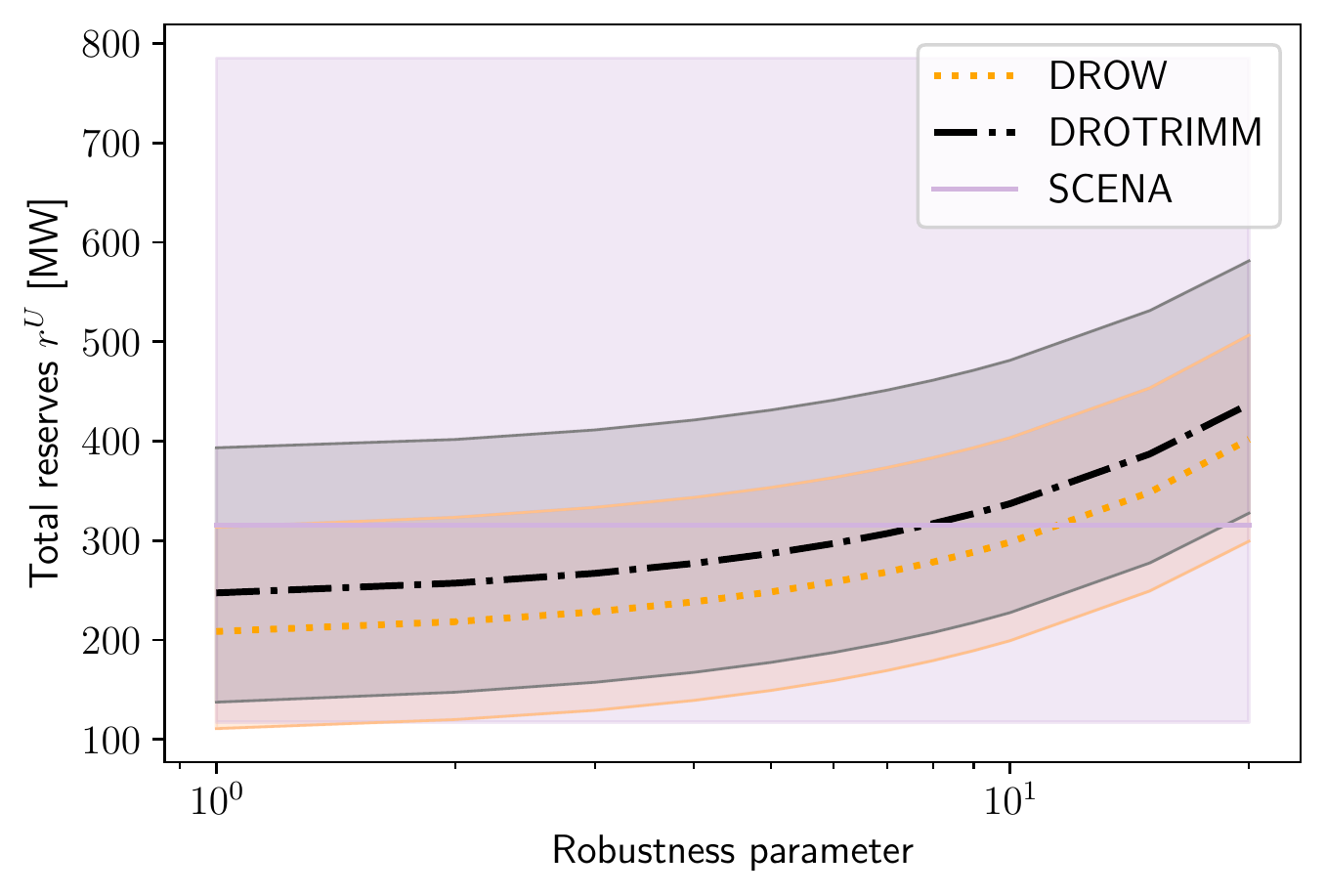}%
  \label{totalreserves_dcopf_sens_100_case_2pos_U}
}

\subfloat[Violation probability (out of sample)]{%
  \includegraphics[width=0.5\textwidth]{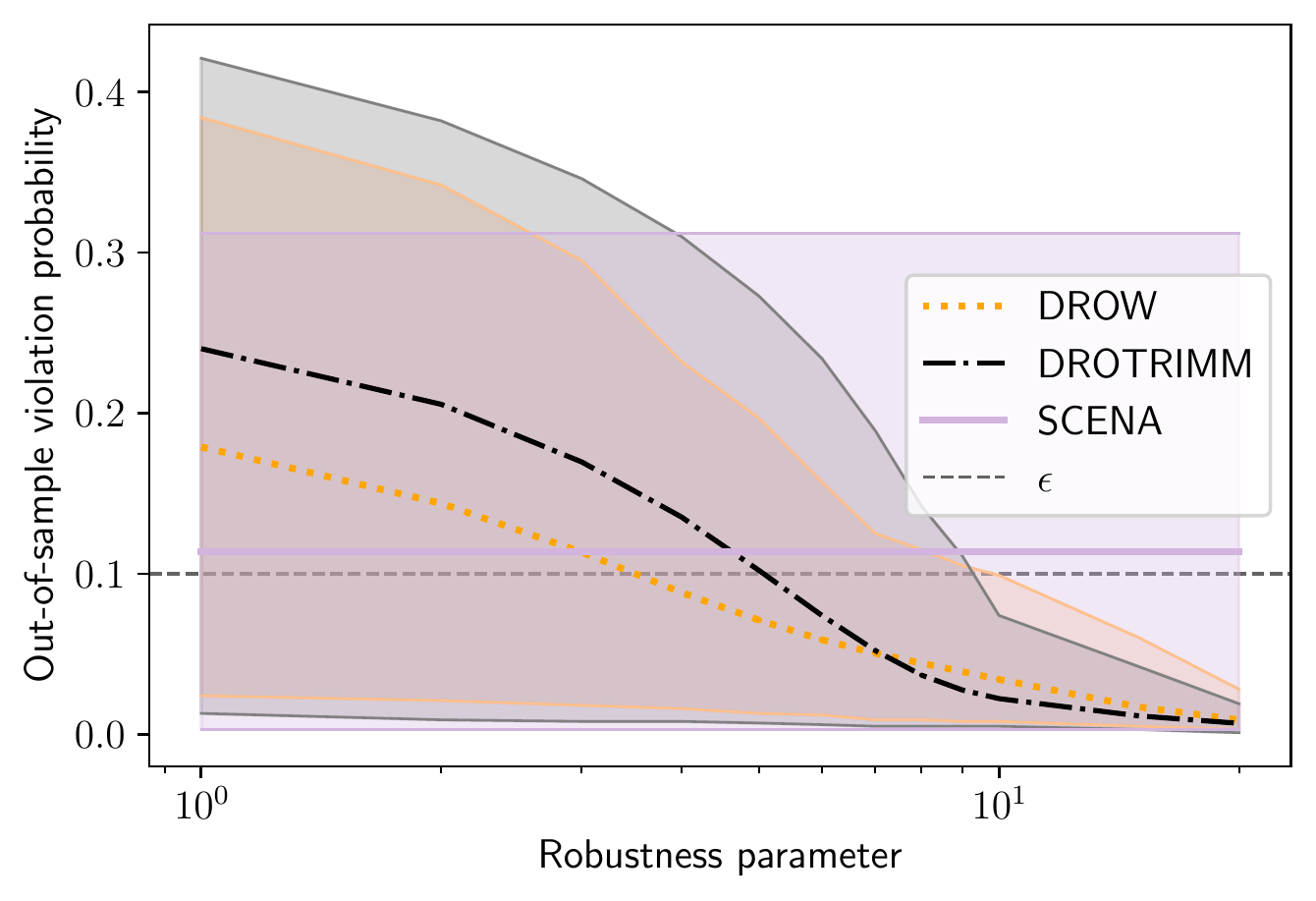}%
  \label{violation_dcopf_sens_100_case2}
}%
\subfloat[Expected system operating cost (out of sample)]{%
  \includegraphics[width=0.5\textwidth]{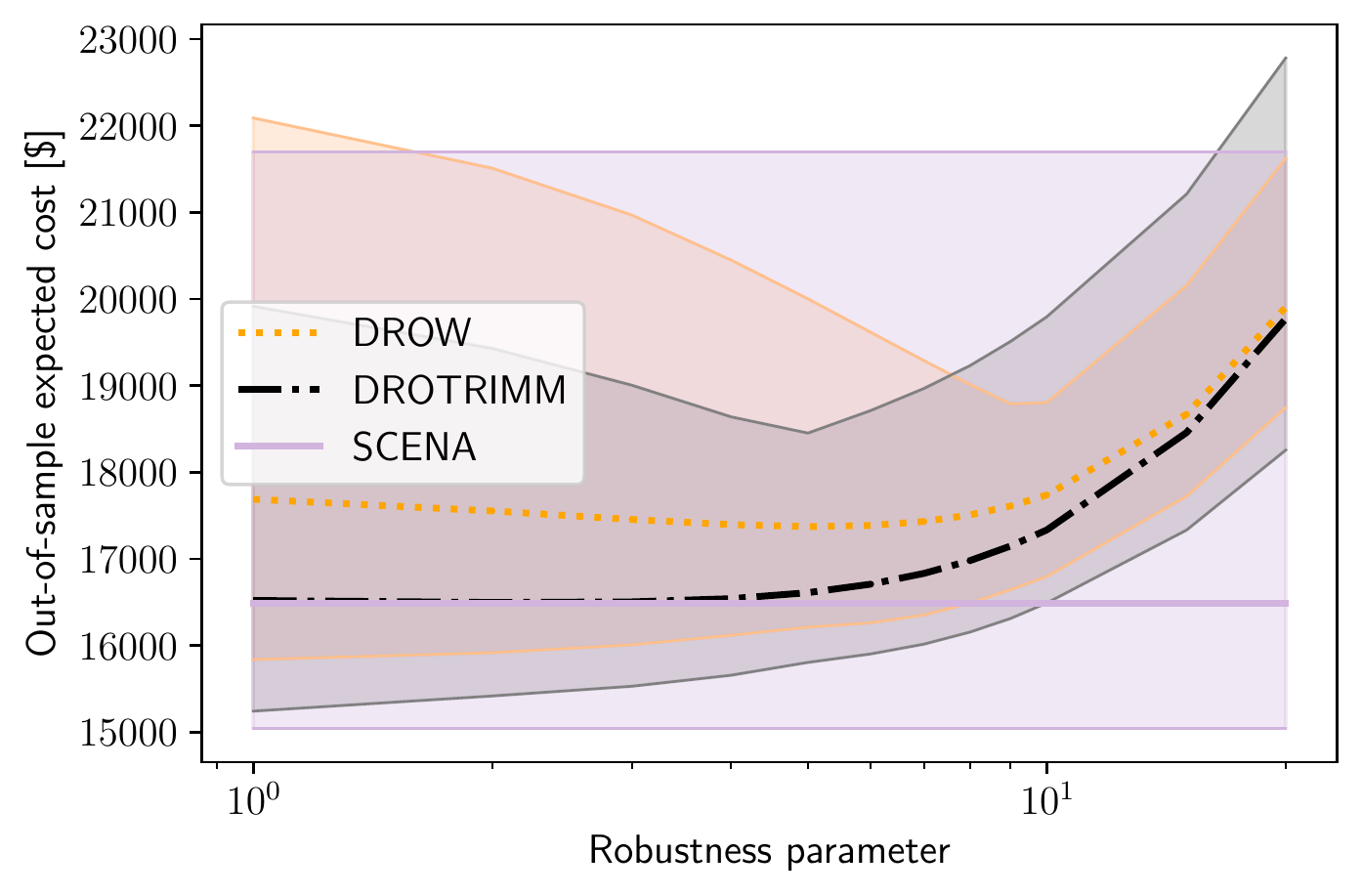}%
  \label{actualexpectedcost_dcopf_sens_100_case2}
}


\caption{High level of wind penetration, $N=100 $ and $\epsilon=0.1$:   Total downward and upward reserve capacity and performance metrics}\label{performance_dcopf_118node_100_case2}

\end{figure}

\begin{figure}
\centering
\subfloat[Total downward reserves]{%
  \includegraphics[width=0.5\textwidth]{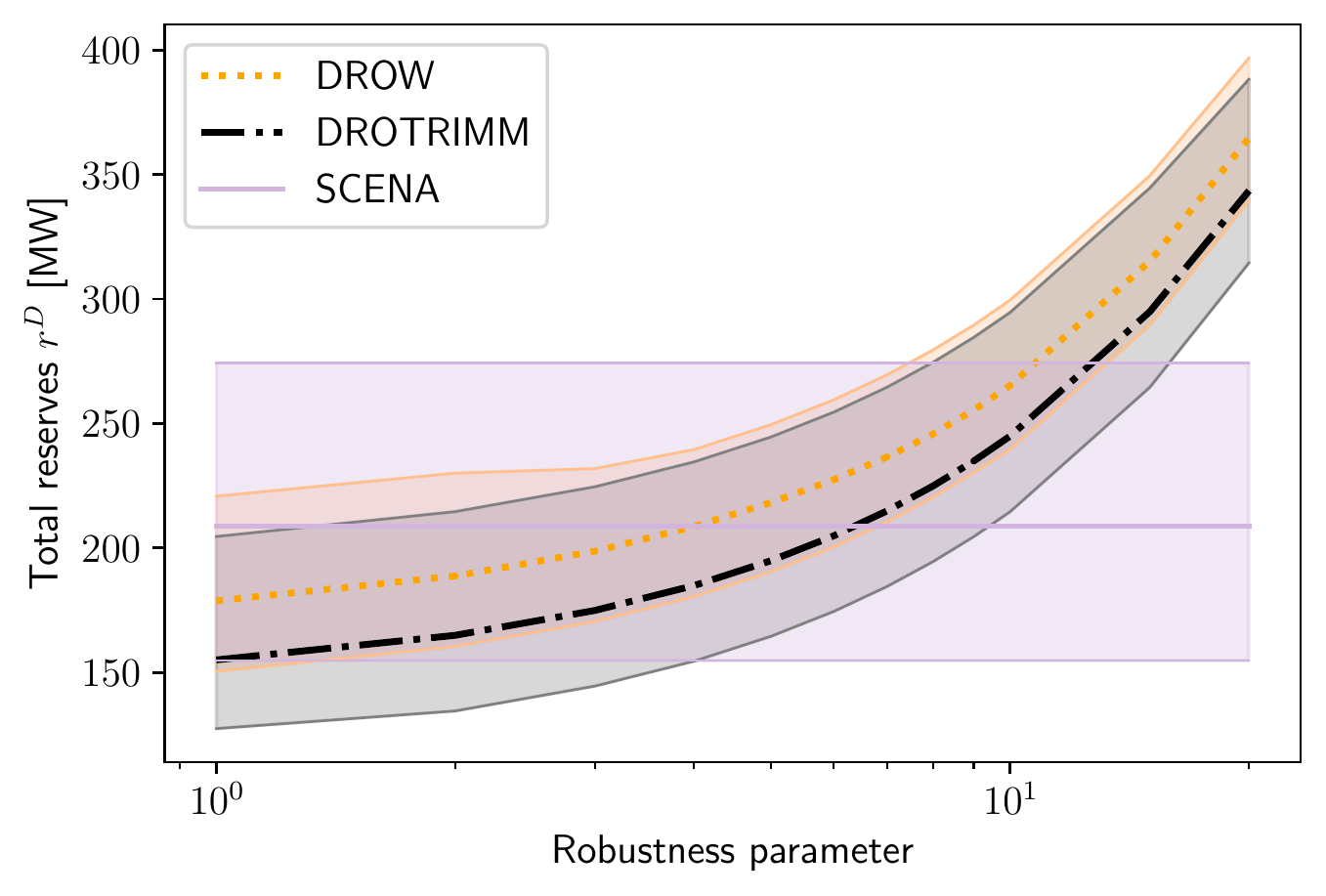}%
  \label{totalreserves_dcopf_sens_300_case_2pos_D}
}%
\subfloat[Total upward reserves]{%
  \includegraphics[width=0.5\textwidth]{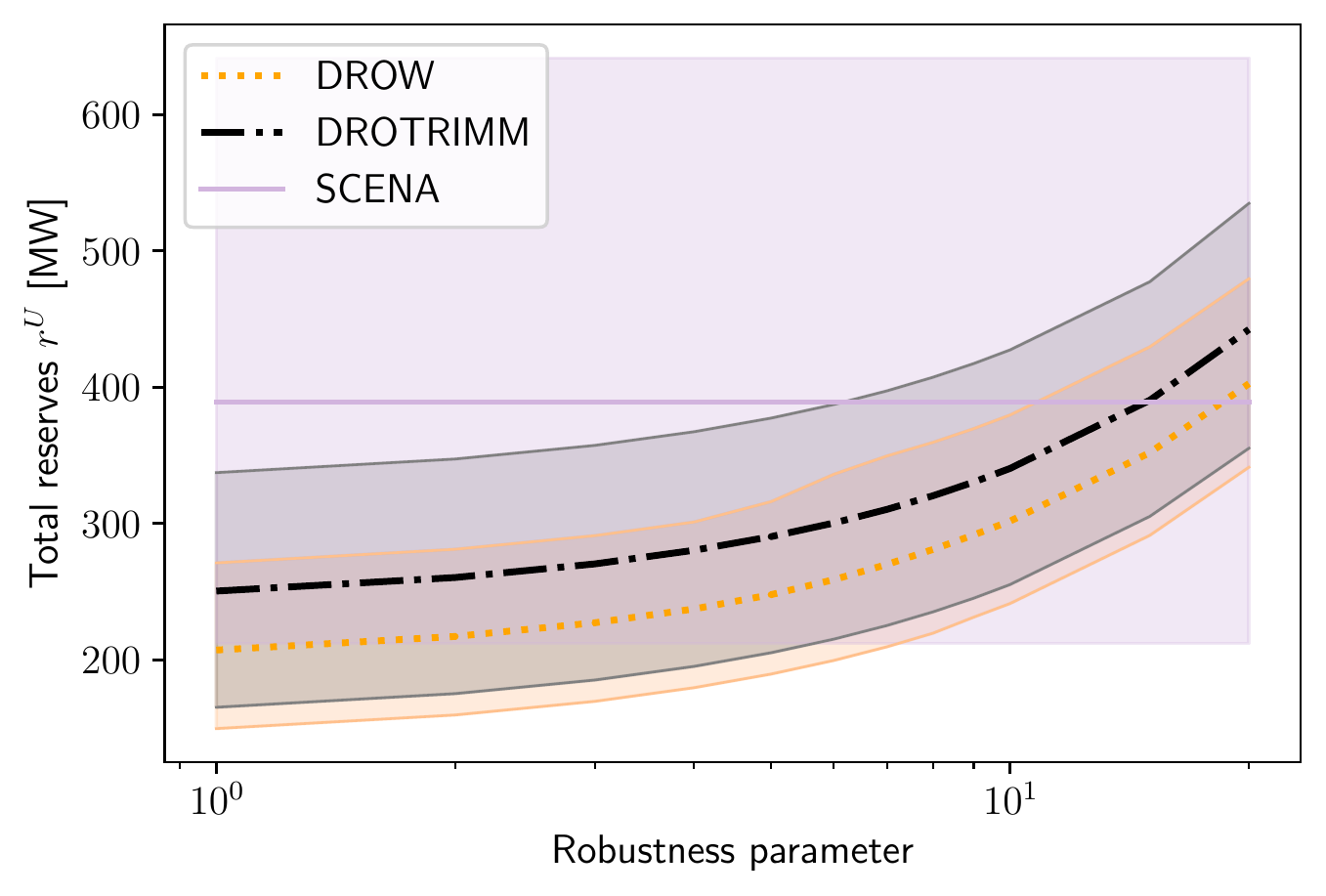}%
  \label{totalreserves_dcopf_sens_300_case_2pos_U}
}

\subfloat[Violation probability (out of sample)]{%
  \includegraphics[width=0.5\textwidth]{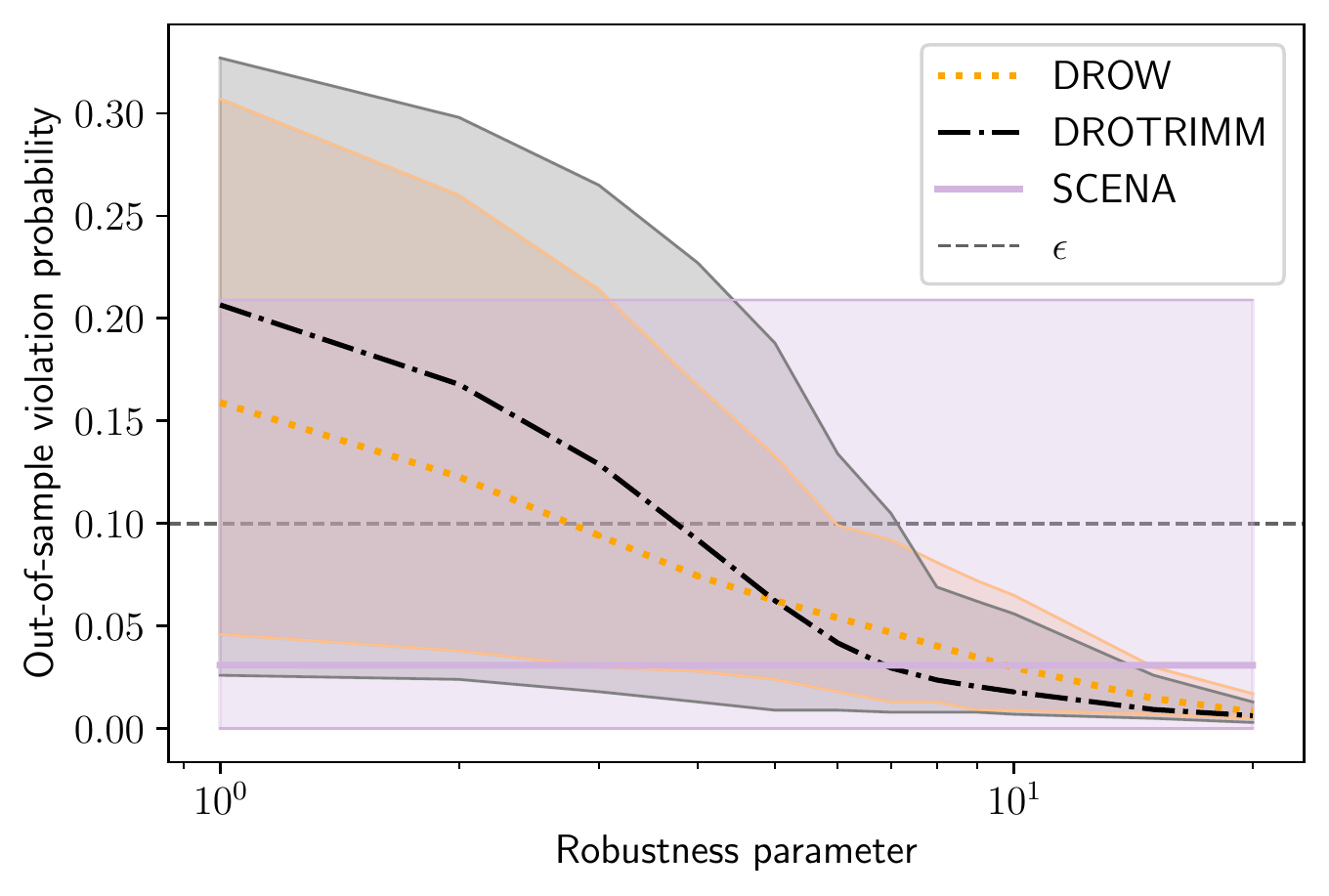}%
  \label{violation_dcopf_sens_300_case2}
}%
\subfloat[Expected system operating cost (out of sample)]{%
  \includegraphics[width=0.5\textwidth]{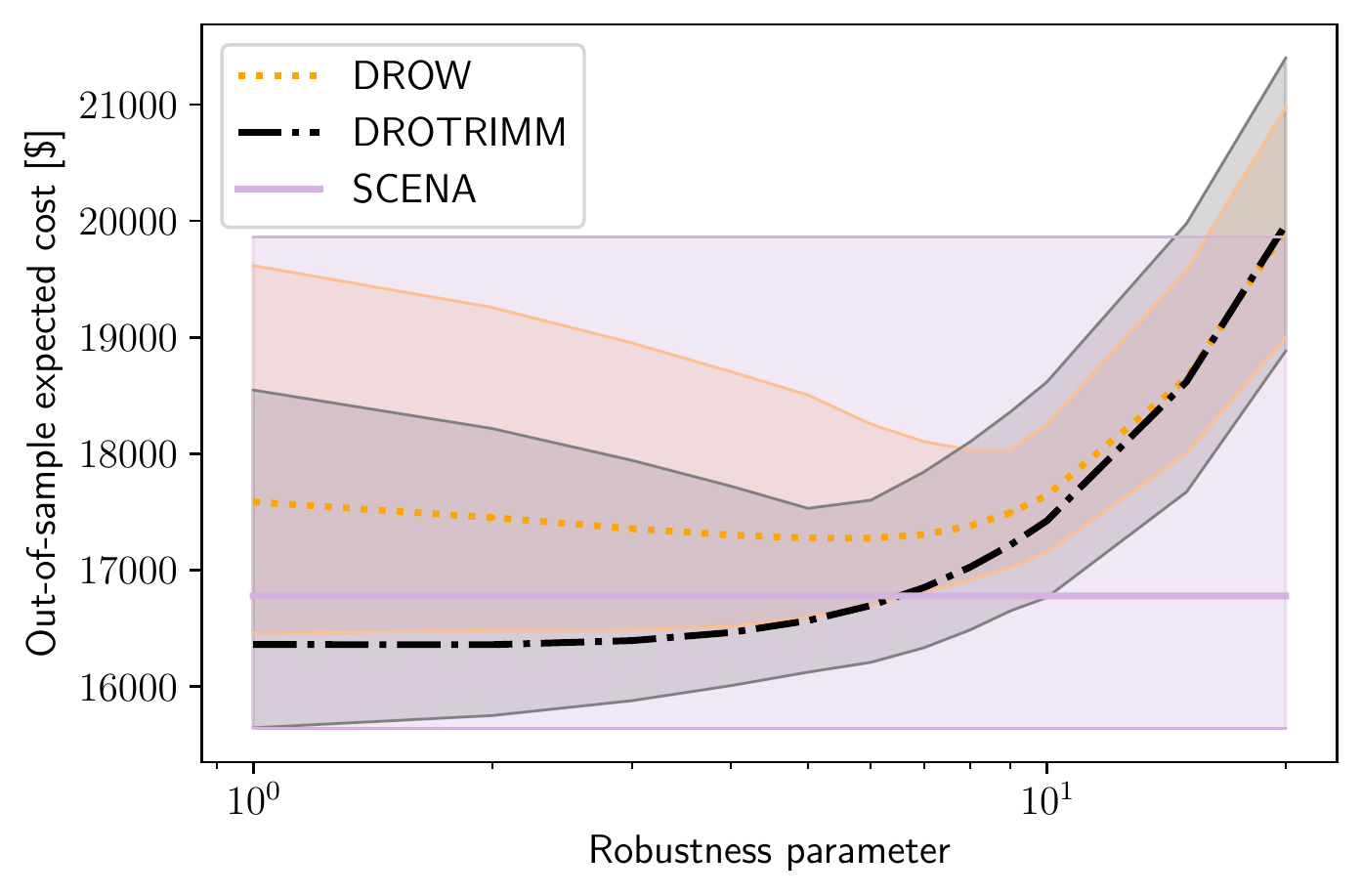}%
  \label{actualexpectedcost_dcopf_sens_300_case2}
}


\caption{High level of wind penetration, $N=300 $ and $\epsilon=0.1$:   Total downward and upward reserve capacity and performance metrics}\label{performance_dcopf_118node_300_case2}

\end{figure}

Figures~\ref{performance_dcopf_118node_100_case2} and \ref{performance_dcopf_118node_300_case2}, and Table~\ref{tab:high} are analogous to Figures~\ref{performance_dcopf_118node_100_case1} and \ref{performance_dcopf_118node_300_case1}, and Table~\ref{tab:medium} of the previous case, respectively. The higher level of wind power penetration in this new instance implies a higher level of uncertainty in the system.  { \color{black} This accentuates the difference in performance between DROW and DROTRIMM when $N = 100$, that is, in a small sample regime. More specifically,  the relative difference between the out-of-sample average expected cost achieved by DROW and DROTRIMM increases from 0.82\% in the previous case to 2.27\% in this new one. It is true, though, that DROW offers reliable OPF solutions with the lowest variance in expected cost when $N=100$, provided that its robustness parameter is optimally tuned, see Table~\ref{tab:high}. However, its superiority in this respect ends when $N$ grows to 300, at which point DROTRIMM provides the most cost-efficient OPF solutions in every respect\footnote{{\color{black}Note in Table~\ref{tab:high} that, while the standard deviation of the expected cost is a bit higher under DROTRIMM than under DROW when $N = 300$, the maximum and minimum values reached by the expected cost under each method reveals that DROTRIMM produces a distribution of the expected cost displaced towards cheaper OPF solutions.}}. Again the reason for this difference in performance has to do with the different provision of upward and downward reserve capacity that  DROW and DROTRIMM prescribe.}

{\color{black} For its part, the SCENA method keeps on providing cheap, but unreliable OPF solutions under a higher level of wind power penetration. In fact, the variability in cost, violation probability and reserves of the OPF solutions given by this method is remarkably high in contrast with that of DROTRIMM and DROW, even higher than in the case of a medium wind power penetration level (compare the range of the box plots in Figure~\ref{performance_dcopf_118node_300_case2})}.


\begin{table}
\caption{High level of wind penetration, summary data for total expected cost [$\$$] under the optimal value of the robustness parameter for DROW and DROTRIMM. }
\label{tab:high}
\centering
\setlength\tabcolsep{3pt} 
{\color{black}
\begin{tabular}{ccccccccc}
\hline
\multirow{2}{*}{} &  &\multicolumn{3}{c}{$N=100$}& &\multicolumn{3}{c}{$N=300$}\\
\cline{3-5} \cline{7-9}
 &    &DROTRIMM &DROW &SCENA &&DROTRIMM &DROW &SCENA \\
\hline
         &  max & 19794  &  18804 &  21699  & &18101 &  18255&  19863 \\
           &avg & 17334      &   17737 &  16483 & & 17025 &17274  & 16779 \\
             &min  & 16490      &  16795&   15040 & & 16486 & 16703  &15639  \\
         &std &   488    & 381  &  1093& & 296&  292 & 767  \\
\hline
\end{tabular}}
\end{table}

 We conclude this section with a remark on computational time. {\color{black} DROTRIMM and DROW have the same complexity (essentially, the number of constraints grows linearly with the sample size $N$)}. The continuous linear program that results from tackling the chance-constrained DRO OPF problem by way of DROTRIMM and the ${\mathbf{CVaR}}$ approximation takes around {\color{black}15 minutes to be solved on average, for a sample size $N=300$, using CPLEX 20.1.0  on a Linux-based server with 22 CPUs   clocking at 2.6 GHz and 200 GB of RAM in total. }

%
\section{Conclusion}\label{sec:conclusion}
In this paper, we have developed a  distributionally robust chance-constrained OPF model that is able to exploit contextual information through an ambiguity set based on probability trimmings. We have provided a reformulation of this model as a continuous linear program using the well-known ${\mathbf{CVaR}}$ approximation. By way of a series of numerical experiments conducted on a modified 118-bus power network with wind uncertainty, we have shown that, by exploiting the statistical dependence between the point forecast of the wind power outputs and its associated forecast error, our approach can identify dispatch solutions that, while satisfying the required system reliability, lead to costs savings of up to several percentage points with respect to the OPF solutions provided by an alternative DRO method that ignores said statistical dependence.

In future work, we plan to devise data-driven schemes for appropriately tuning the robustness parameter in our distributionally robust chance-constrained OPF model in accordance with the risk preferences of the system operator (for instance, by resorting to cross-validation or bootstrapping). We also want to extend this model to account for intertemporal constraints, which, among other things, would involve adapting our probability-trimming-based ambiguity set to deal with stochastic processes and time series data.



\section*{Acknowledgments}

This work was supported in part by the European Research Council (ERC) under the EU Horizon 2020 research and innovation program (grant agreement No. 755705), in part by the Spanish Ministry of Science and Innovation through project PID2020-115460GB-I00/AEI/10.13039/501100011033, and in part by the Junta de Andalucía (JA) and the European Regional Development Fund (FEDER) through the research project P20\_00153. Finally, the authors thankfully acknowledge the computer resources, technical expertise, and assistance provided by the SCBI (Supercomputing and Bioinformatics) center of the University of M\'alaga.

\bibliography{References.bib}

\newpage

\appendix
\section{Notation}\label{appendix_notation}
The main notation used throughout the text is stated below for quick reference.  Other symbols are defined as required.
{
\color{black}

\subsection{Sets, numbers and indices}
\begin{ldescription}{$xxx$}
\item [${\color{black}\mathcal{B}}$] Set of buses, indexed by $b$.
\item [$\mathcal{L}$] Set of lines, indexed by $\ell$.
\item [$\mathcal{G}$] Set of generators (dispatchable units), indexed by $j$.
\item [$\mathcal{W}$] Set of wind power plants, indexed by $m$.
\end{ldescription}
\subsection{Parameters and functions}
\begin{ldescription}{$xxxxxxx$}
\item [$\mathbf{f}$] Array of forecasted power outputs [MW].
\item [$\widetilde{\mathbf{f}}$] Array of nominal (p.u.) forecasted power outputs.
\item [$\mathbf{L}$] Array of loads [MW].
\item [$\mathbf{g}^{\min}, \mathbf{g}^{\max}$] Array of upper and lower capacity limits of generators [MW].
\item [$\textbf{Cap}$] Array of line capacities [MW].
\item [$\overline{\textbf{C}}$] Array of installed capacities of the wind power plants [MW].

\item [${\color{black}\mathbf{M}^{\mathcal{G}/\mathcal{W}/\mathcal{B}}}$ ] Matrix of DC power transfer distribution factors, which maps nodal power injections to line flows for generators/wind farms/loads.
\item [$\mathbf{c}^{D}, \mathbf{c}^{U}$ ] Array of downward and upward reserve capacity costs [\$/MW].
\item [$C(\cdot)$] Total production cost function,  which is  given
    by the sum of $|\mathcal{G}|$ convex piecewise linear cost functions  with $S_j$ pieces for generator $j$, i.e.,
  $C(\tilde{\mathbf{g}}(\boldsymbol{\omega})):=\sum_{j \in \mathcal{G}} \max_{s=1,\ldots, S_j } \{ m_{js}\tilde{g}_j(\boldsymbol{\omega}) +n_{js} \}$,
    where $m_{js}, n_{js}$ are  the slope and the intercept of the $s$-th piece   for generator $j$, respectively [\$].

\end{ldescription}

\subsection{Random variables and uncertain parameters}
\begin{ldescription}{$xxxxxxl$}
\item [$\mathbf{z}$] Random vector of features/covariates.
\item [$\boldsymbol{\omega}$] Random vector  representing the wind power forecast errors of the ${\color{black} |\mathcal{W}|}$ wind power plants [MW].

\item [$\boldsymbol{\xi}$] Random vector representing the pair of features/covariates and the wind power forecast errors of the ${\color{black} |\mathcal{W}|}$ wind farms, that is, $\boldsymbol{\xi}:=(\mathbf{z},\boldsymbol{\omega})$.
\item [$W_m$] Actual power output at wind power plant $m\in \mathcal{W}$ in per unit.
\item [$\Xi_{\omega}$] Support set of the random vector $\boldsymbol{\omega}$.
\item [$\widetilde{\Xi}_{\boldsymbol{\omega}}$] Support set of the random vector $\boldsymbol{\omega}$ conditional on $\mathbf{z} = \mathbf{f}$, which is given by the hypercube $\prod_{m \in \mathcal{W}}  [-f_m,\overline{C}_m-f_m]$.

\item [$\Xi$] Support set of the random vector $(\mathbf{z},\boldsymbol{\omega})$.
\item [$\widetilde{\Xi}$] Contextual information, that is, the event $( \mathbf{z}=\mathbf{f};\  \boldsymbol{\omega}\in \widetilde{\Xi}_{\boldsymbol{\omega}})$.
\item [$\Omega$]Random variable defined as $\sum_{{\color{black}m\in \mathcal{W}}}\omega_m$,   which describes the system-wise aggregate wind power forecast error [MW].
\item [$\widetilde{\Xi}_{\Omega}$] Contextual information linked to the random vector $(\mathbf{z},\Omega)$, that is, the event $ (\mathbf{z}=\mathbf{f}; \Omega \in [\underline{\Omega},\overline{\Omega}]), \; \text{with}\; [\underline{\Omega},\overline{\Omega}]=\left[-\sum_{m\in \mathcal{W}}f_m,\sum_{m\in \mathcal{W}}(\overline{C}_m-f_m)\right]$.
\item [$\tilde{\mathbf{g}}(\boldsymbol{\omega})$] Array of power generation outputs of generators (random vector) [MW].
\item [$\tilde{\mathbf{r}}(\boldsymbol{\omega})$] Array of reserves deployed by generators (random vector) [MW].
\item [$\mathbb{E}_{Q}$] Expectation operator with respect to the probability measure $Q$.
\item [$\delta_{\xi}$] Dirac distribution at $\xi$.
\end{ldescription}

%

%
\subsection{Variables}
\begin{ldescription}{$xxxxxxl$}
\item [$\mathbf{g}$] Generators' power dispatch [MW].
\item [$\boldsymbol{\beta}$] Array of generators' participation factors.
\item [$\mathbf{r}^{D},\mathbf{r}^{U}$] Array of downward/upward  reserve  capacities  provided by generators  [MW].
\item [$\mathbf{x}$] Vector of  decision variables, that is, $\mathbf{x}:=(\mathbf{g},\boldsymbol{\beta},\mathbf{r}^D,\mathbf{r}^U)$.
\item [$\mathbf{y}$] Vector of  first-stage decision variables (power dispatch and reserve capacity provision), that is, $\mathbf{y}:=(\mathbf{g},\mathbf{r}^D,\mathbf{r}^U)$.
\end{ldescription}

\subsection{Other symbols}\label{appendix:othersymbols}
\begin{ldescription}{$xxxxxxl$}
\item [$\mathbf{1}$] Array of ones (of appropriate dimension).
\item [$\mathbf{0}$] Array of zeros (of appropriate dimension).
\item [$|A|$] Cardinality of a set $A$.
\item [$(x)^+$] Positive part of $x$, i.e.,  $\max\{x,0 \}$.
\item [$\lfloor x \rfloor$] Floor function of $x$, given by $\max\{ m \in \mathbb{Z} \;/\; m\leqslant x\}$.
\item [$\langle \cdot, \cdot\rangle$ ] Dot product.
\item [$W_1$] $1$-Wasserstein distance.
\item [$\mathcal{P}_1(\Xi),\mathcal{P}_1(\widetilde{\Xi})$] The set of all probability distributions  with finite first moment supported on $\Xi,\widetilde{\Xi}$, respectively.
\item [$\mathcal{R}_{1-\alpha}(P)$] The set of
all $(1-\alpha)$-trimmings  of the probability distribution $P$.
\item [$\rho$] Robustness parameter.
\item [$\underline{\rho}_{N\alpha}$] Minimum transportation budget.
\item [$S_B$] Support function of a set $B \subseteq \mathbb{R}^d$, defined as $S_B(a):=\sup_{b \in B }\langle a,b\rangle$.
\item [$Q-{\mathbf{CVaR}}_{\epsilon}(\phi(\boldsymbol{\omega}))$] Conditional Value at Risk at level $\epsilon \in (0,1)$ of $\phi(\boldsymbol{\omega})$
under the probability measure $Q$; that is,   the value  $\inf_{\tau \in \mathbb{R}}\{\tau +\frac{1}{\epsilon}\mathbb{E}_{Q}[(\phi(\boldsymbol{\omega})-\tau)^{+}] \}$.

\end{ldescription}

}

\section{Real-time re-dispatch problem}\label{appendix_real_time_formulation}
This appendix contains the optimization program used to evaluate the out-of-sample performance of a given solution of the chance-constrained DC-OPF problem. Given $N$, a data-driven solution $\mathbf{y}_N:=(\mathbf{g},\mathbf{r}^D,\mathbf{r}^U)_N$ and  a realization of the forecast error $\widehat{\boldsymbol{\omega}}_i$, the operator of the system solves the following deterministic linear program:
{\color{black}
\begin{align}
\min_{\mathbf{r},\boldsymbol{\Delta d},\boldsymbol{\Delta \omega}}&     C(\mathbf{g}_N+\mathbf{r})+\langle \mathbf{c}^{\text{shed}},\boldsymbol{\Delta d}  \rangle+\langle \mathbf{c}^D, \mathbf{r}^D_N \rangle +\langle \mathbf{c}^U, \mathbf{r}^U_N \rangle \label{obj_real_dispatch}\\
{\text {s.t.}}\,&\  \mathbf{0}\leqslant \boldsymbol{\Delta d} \leqslant \mathbf{L} \label{real_dispatch_LS}\\
&\mathbf{0}\leqslant \boldsymbol{\Delta \omega} \leqslant \mathbf{f} +\widehat{\boldsymbol{\omega}}_i \label{real_dispatch_WS}\\
&-\mathbf{r}^D_{N} \leqslant \mathbf{r}\leqslant \mathbf{r}^U_{N} \label{real_dispatch_DR}\\
&\langle\mathbf{1},\mathbf{r}\rangle+\langle\mathbf{1},\boldsymbol{\Delta d}\rangle+\langle\mathbf{1},\widehat{\boldsymbol{\omega}}_i-\boldsymbol{\Delta \omega}\rangle=0 \label{real_dispatch_PB}\\
&-\textbf{Cap}\leqslant  \mathbf{M}^{\mathcal{G}}(\mathbf{g}_N+\mathbf{r})+\mathbf{M}^{\mathcal{W}}(\mathbf{f}+\widehat{\boldsymbol{\omega}}_i-\boldsymbol{\Delta \omega})-\mathbf{M}^{\mathcal{B}}(\mathbf{L}-\boldsymbol{\Delta d})  \leqslant \textbf{Cap} \label{real_dispatch_TC}
\end{align}}
where $\mathbf{r}$, $\boldsymbol{\Delta d}$ and  $\boldsymbol{\Delta \omega}$ are the deployed reserves, load shedding and wind spillage vector of decision variables; and the parameter $c_{b}^{\text{shed}} $ is the load shedding cost at bus $b$.
The objective function in \eqref{obj_real_dispatch} minimizes the total operational cost of the system, which comprises the electricity generation cost, the load shedding cost and the total cost of up- and down-reserve capacities. The latter is known and constant and thus, does not intervene in the minimization. Constraints~\eqref{real_dispatch_LS} and~\eqref{real_dispatch_WS} limit the amount of load involuntarily curtailed and the amount of wind power unused to the actual realizations of the load and the wind power production, respectively. Constraint~\eqref{real_dispatch_DR} ensures that the deployed reserves are kept within the reserve capacities scheduled in the forward stage. Constraint~\eqref{real_dispatch_PB} constitutes the real-time power balance equation and, finally, constraints~\eqref{real_dispatch_TC} enforce the transmission capacity limits.



\section*{Supplementary material}

## 1. Proof of Proposition 1

The proof of this proposition is a direct application of [8, Theorem 1], after noticing that the inner supremum in constraint (14) involves a maximum of linear functions. Variables $\mathbf{v}_{ik}$ and $\boldsymbol{\gamma}_{ik}$ in optimization problem (16) are auxiliary variables that result from the dualization of the inner supremum that appears in the **CVaR** approximation of the chance constraint system, that is, in (13). This dualization is critical to the reformulation of (13) as a tractable mathematical program, see [8, Theorem 1].

## 2. Proof of Proposition 2

Based on [8, Theorem 1], the DRO problem defined by (19) can be reformulated as follows:

$$\inf_{\lambda \geqslant 0; \overline{\mu}_i, \forall i \leqslant N; \theta \in \mathbb{R}} \quad \lambda \rho + \theta + \frac{1}{N\alpha} \sum_{i=1}^{N} \overline{\mu}_i + \langle \mathbf{c}^D, \mathbf{r}^D \rangle + \langle \mathbf{c}^U, \mathbf{r}^U \rangle \tag{SM-1a}$$

$$\text{s.t.} \ \overline{\mu}_i + \theta + \lambda \|\mathbf{z}^* - \widehat{\mathbf{z}}_i\| \geqslant \sup_{\Omega \in [\underline{\Omega}, \overline{\Omega}]} \left( \sum_{j \in \mathcal{G}} \max_{s=1,\dots,S_j} \left\{ m_{js} \left[ g_j - \beta_j \Omega \right] + n_{js} \right\} \right.$$
$$\left. -\lambda |\Omega - \widehat{\Omega}_i| \right), \ \forall i \leqslant N \tag{SM-1b}$$

$$\overline{\mu}_i \geqslant 0, \ \forall i \leqslant N \tag{SM-1c}$$

Now, constraint (SM-1b) is equivalent to the following ones:

$$\overline{\mu}_i + \theta + \lambda \|\mathbf{z}^* - \widehat{\mathbf{z}}_i\| \geqslant t_i, \ \forall i \leqslant N \tag{SM-2a}$$

$$t_i \geqslant \sup_{\Omega \in [\underline{\Omega}, \overline{\Omega}]} \sum_{j \in \mathcal{G}} \max_{s=1,\dots,S_j} \left\{ m_{js} \left[ g_j - \beta_j \Omega \right] + n_{js} \right\} - \lambda |\Omega - \widehat{\Omega}_i|, \ \forall i \leqslant N \tag{SM-2b}$$

In order to reformulate the supremum on the right-hand side of (SM-2b), we resort to the following partition of the set $\{1, \dots N\}$:

$$\underline{I} := \{i \in \{1, \dots N\} : \widehat{\Omega}_i < \underline{\Omega}\} \tag{SM-3a}$$

$$I := \{i \in \{1, \dots N\} : \widehat{\Omega}_i \in [\underline{\Omega}, \overline{\Omega}]\} \tag{SM-3b}$$

$$\overline{I} := \{i \in \{1, \dots N\} : \widehat{\Omega}_i > \overline{\Omega}\} \tag{SM-3c}$$

In this way, because of the convexity of the sum of a maximum of affine functions, constraint (SM-2b) can be replaced by the following set of constraints:

$$t_i \geqslant \sum_{j \in \mathcal{G}} \max_{s=1,\dots,S_j} \left\{ m_{js} \left[ g_j - \beta_j \underline{\Omega} \right] + n_{js} \right\} - \lambda(\underline{\Omega} - \widehat{\Omega}_i), \ \forall i \in \underline{I} \tag{SM-4a}$$

$$t_i \geqslant \sum_{j \in \mathcal{G}} \max_{s=1,\dots,S_j} \left\{ m_{js} \left[ g_j - \beta_j \overline{\Omega} \right] + n_{js} \right\} - \lambda(\overline{\Omega} - \widehat{\Omega}_i), \ \forall i \in \underline{I} \tag{SM-4b}$$

$$t_i \geqslant \sum_{j \in \mathcal{G}} \max_{s=1,\dots,S_j} \left\{ m_{js} \left[ g_j - \beta_j \underline{\Omega} \right] + n_{js} \right\} + \lambda(\underline{\Omega} - \widehat{\Omega}_i), \ \forall i \in \overline{I} \tag{SM-4c}$$

$$t_i \geqslant \sum_{j \in \mathcal{G}} \max_{s=1,\dots,S_j} \left\{ m_{js} \left[ g_j - \beta_j \overline{\Omega} \right] + n_{js} \right\} + \lambda(\overline{\Omega} - \widehat{\Omega}_i), \ \forall i \in \overline{I} \tag{SM-4d}$$

$$t_i \geqslant \sum_{j \in \mathcal{G}} \max_{s=1,\dots,S_j} \left\{ m_{js} \left[ g_j - \beta_j \overline{\Omega} \right] + n_{js} \right\} - \lambda(\overline{\Omega} - \widehat{\Omega}_i), \ \forall i \in I \tag{SM-4e}$$

$$t_i \geqslant \sum_{j \in \mathcal{G}} \max_{s=1,\dots,S_j} \left\{ m_{js} \left[ g_j - \beta_j \underline{\Omega} \right] + n_{js} \right\} + \lambda(\underline{\Omega} - \widehat{\Omega}_i), \ \forall i \in I \tag{SM-4f}$$

$$t_i \geqslant \sum_{j \in \mathcal{G}} \max_{s=1,\dots,S_j} \left\{ m_{js} \left[ g_j - \beta_j \widehat{\Omega}_i \right] + n_{js} \right\}, \ \forall i \in I \tag{SM-4g}$$

Introducing epigraphical auxiliary variables $\underline{t}_{ij}, \overline{t}_{ij}$ and $\widehat{t}_{ij}$, we finish the proof. $\quad\square$.

## 3. Illustrative example (3-bus system)

In this appendix, we use a small three-node system to illustrate how our DRO framework based on probability trimmings, named DROTRIMM, compares to that proposed in [3], which we refer to as KNNDRO. This other approach is based on a local inference method (specifically, a $K$-nearest neighbors), to construct, from the

joint data sample, the conditional empirical distribution at which a Wasserstein ball is centered. We also compare DROTRIMM with the DROW approach introduced in Section 5 of the main text. Actually, DROW becomes equivalent to KNNDRO when taking $K = N$.

The topology of the three-bus system has been taken from [20]. It includes three lines connecting buses 1–2, 2–3, and 1–3, three generators located at nodes 1, 2 and 3, and a 200-MW load connected to bus 3. The production cost of the three generators is modeled as a piecewise function consisting of three pieces. Further details on the generators' and network's parameters can be found in 4.

We consider one wind power plant placed at bus 2 with an installed capacity of $\overline{C}_1 = 60$ MW. Its predicted power output is assumed to be $z^* = f_1 = 30$ MW in this example. We also assume that the system operator has a series of data pairs given by past point forecasts of the power output at the wind farm and their associated forecast error. Figure 1 shows a heat map of the true bivariate joint distribution of the forecast power output and its error, together with a kernel estimate of the probability density function of the random forecast error $\omega$ conditional on $z^* = 30$ MW. The joint support set $\Xi$ is polyhedral, and recall that the conditional support varies with the value of $z$ (see the red line in Figure 1).

The box plots corresponding to the total downward and upward reserve capacity that is procured, the violation probability and the expected cost delivered out of sample by each of the considered CC-DRO OPF models is depicted in Figures 2 and 4 as a function of their respective robustness parameter, estimated over 200 independent runs for a fixed sample size $N = 30$ and $N = 2000$, respectively. Similarly, the plots in Figures 3 and 5 pertain to the generators' dispatch and their participation factors. The robustness parameter of KNNDRO and DROW corresponds to the radius of the Wasserstein ball these methods use as the ambiguity set. For its part, the

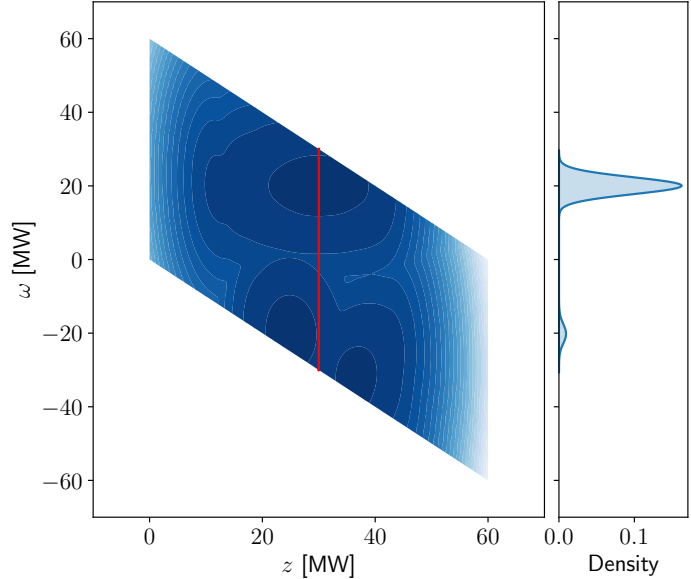

Figure 1. Heat map of the true joint distribution and kernel estimate of the true conditional density given $z^* = 30$ MW

robustness parameter for DROTRIMM is to be greater than or equal to the minimum transportation budget (see Definition 4). Therefore, what the system operator really needs to tune in DROTRIMM is the budget excess as done in [8]. This is what we represent in the $x$-axes of the aforementioned figures for this method.

As in the case study in Section 5 of the main text, the color-shaded areas have been obtained by joining the minimum and maximum edge cases of the box plots, while the associated bold colored lines link their means. The number of neighbors $K_N$ we have considered for KNNDRO is given by the logarithmic rule, that is, if $N$ is the sample size of the joint data, then the number of neighbors is computed as $K_N = \lfloor N/(\log(N+1)) \rfloor$. To ensure a fair comparison, we have also taken $\alpha_N = K_N/N$ for DROTRIMM. Figures 2–5 provide information on the ability of each method to discover good dispatch solutions, that is, scheduling plans for power production and

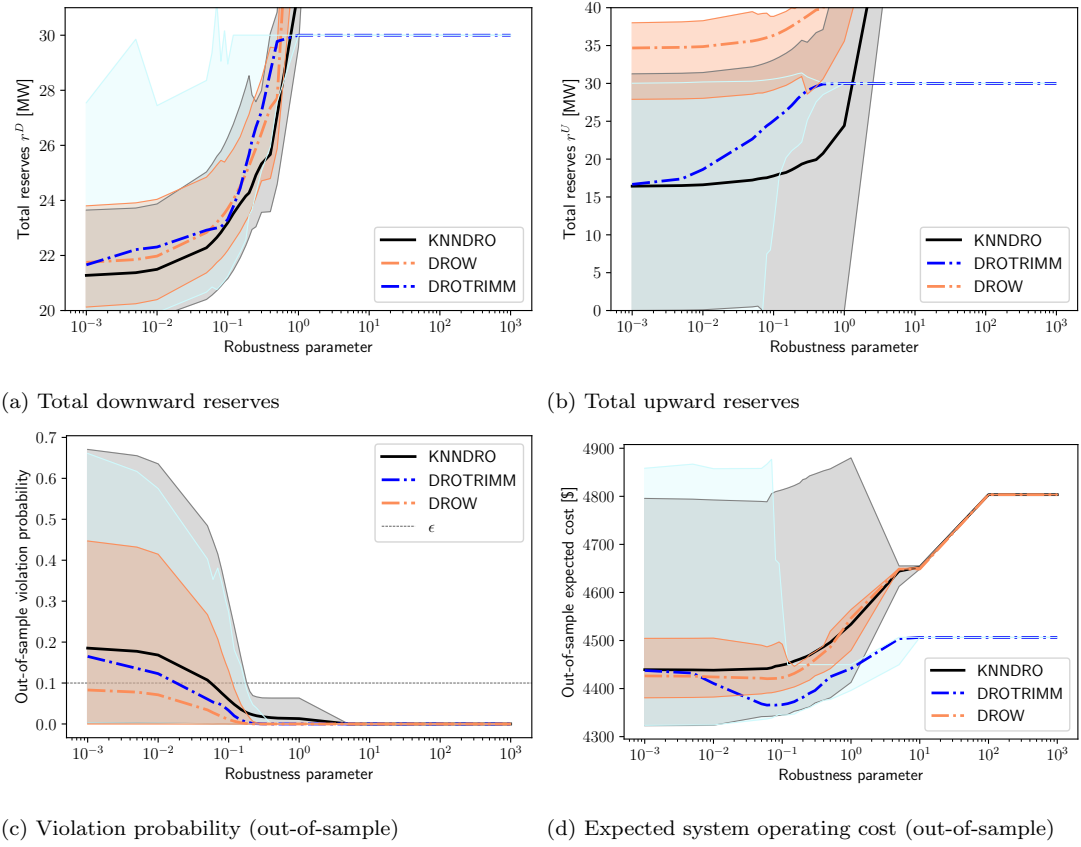

(a) Total downward reserves

(b) Total upward reserves

(c) Violation probability (out-of-sample)

(d) Expected system operating cost (out-of-sample)

Figure 2. Three-bus system, sample size $N = 30$ and $\epsilon = 0.1$: Total downward and upward reserves and performance metrics

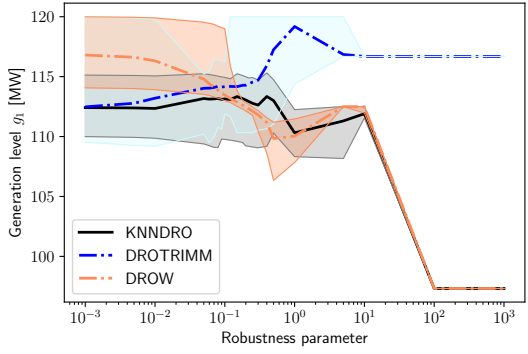

(a) $g_1$

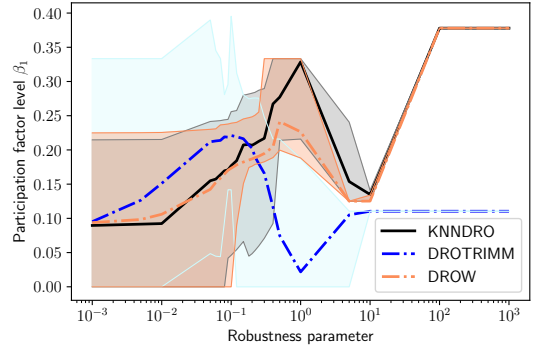

(b) $\beta_1$

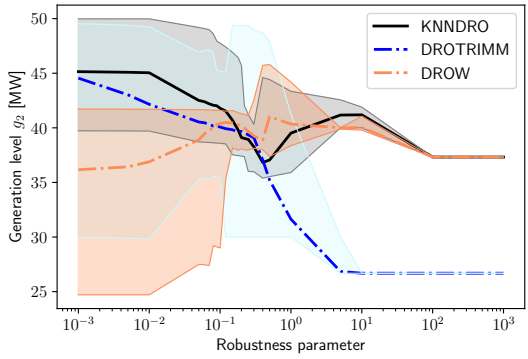

(c) $g_2$

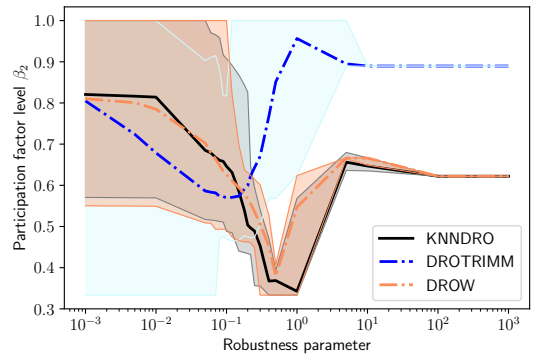

(d) $\beta_2$

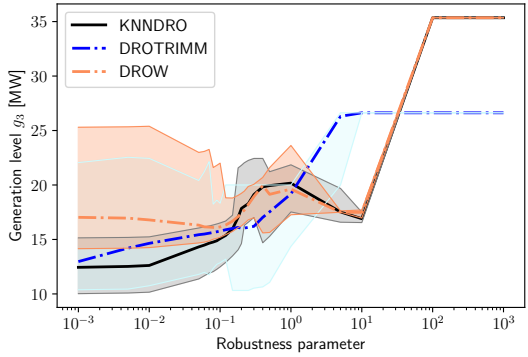

(e) $g_3$

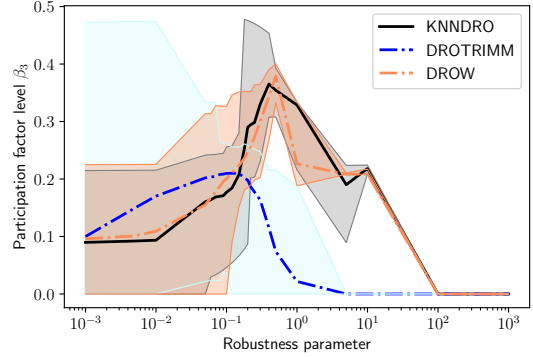

(f) $\beta_3$

Figure 3. Three-bus system, sample size $N = 30$ and $\epsilon = 0.1$: Generators' dispatch and participation factors

6

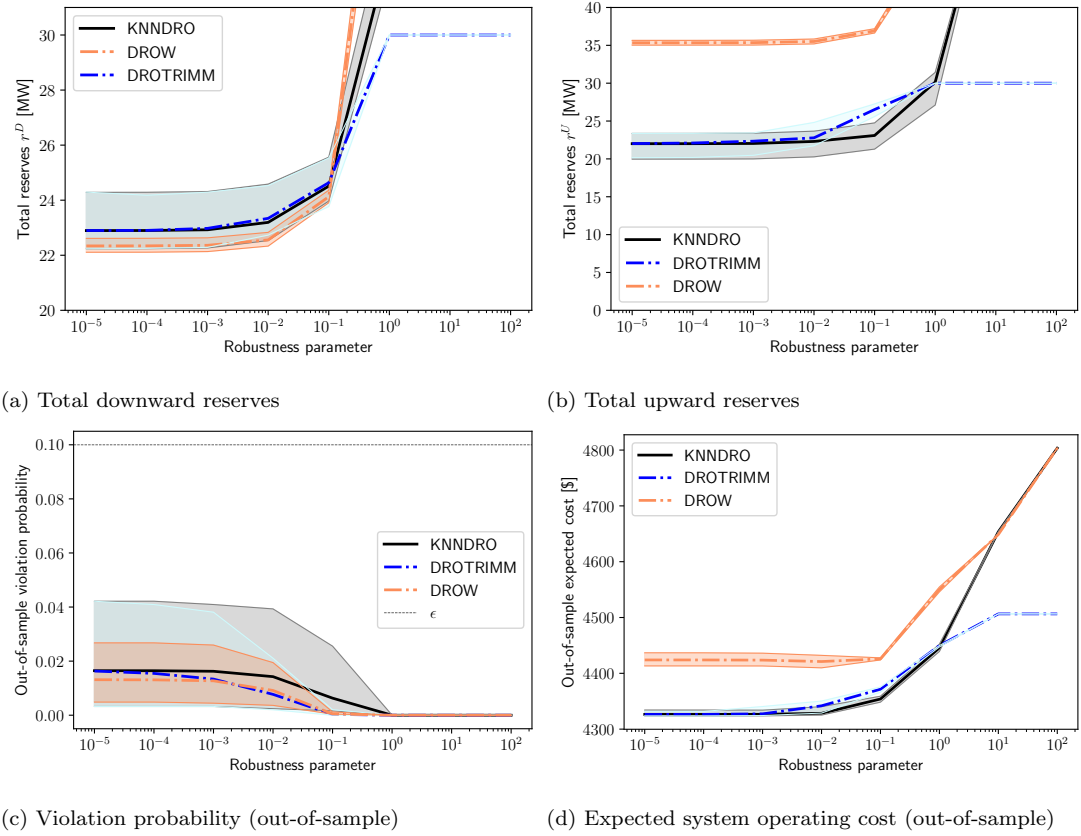

(a) Total downward reserves

(b) Total upward reserves

(c) Violation probability (out-of-sample)

(d) Expected system operating cost (out-of-sample)

Figure 4. Three-bus system, sample size $N = 2000$ and $\epsilon = 0.1$: Total downward and upward reserves and performance metrics

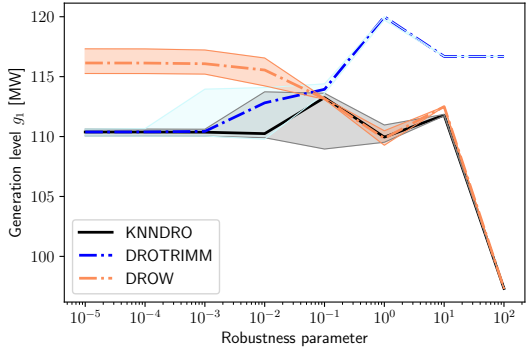

(a) $g_1$

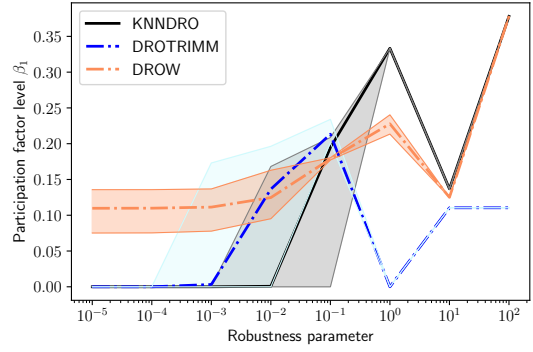

(b) $\beta_1$

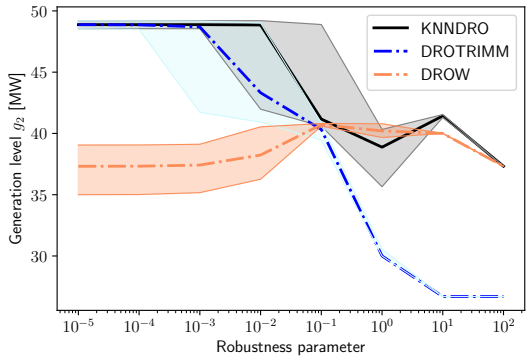

(c) $g_2$

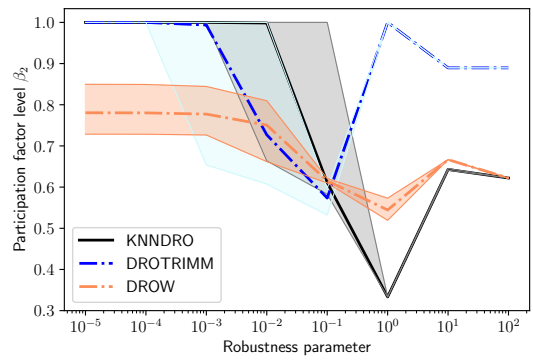

(d) $\beta_2$

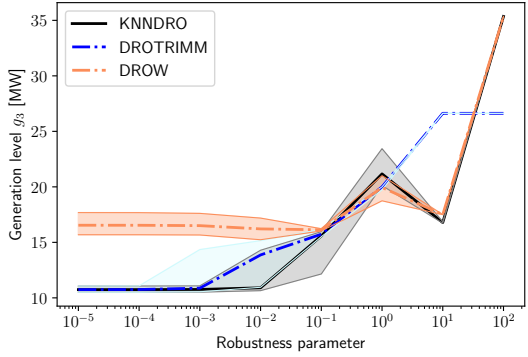

(e) $g_3$

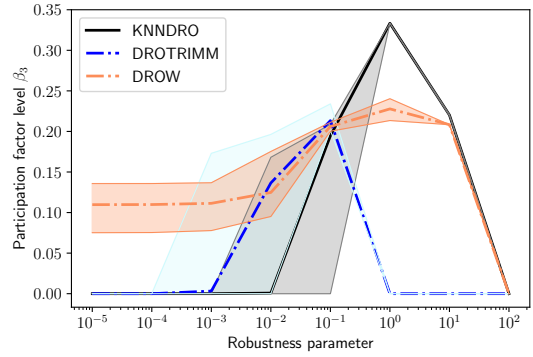

(f) $\beta_3$

Figure 5. Three-bus system, sample size $N = 2000$ and $\epsilon = 0.1$: Generators' dispatch and participation factors

8

reserve capacity provision that are cost-efficient in expectation while guaranteeing the desired system reliability. We consider two settings, namely, the case of a small sample size ($N = 30$), for which the available joint data is expected to carry little information on the true conditional distribution of the wind power forecast error (Figures 2 and 3), and another one where the sample size is notably higher, specifically $N = 2000$ (Figures 4 and 5).

Figure 2 shows that, for a value of the robustness parameter equal to or greater than $10^2$, DROTRIMM recovers the *robust* dispatch, that is, the dispatch that performs the best under the worst-case value of the forecast error, having predicted 30 MW of wind power production. On the contrary, KNNDRO and DROW recover the *unconditional robust* dispatch needing a substantially higher amount of reserves than the ones required by DROTRIMM. This is so because the probability distributions in the Wasserstein balls that DROW and KNNDRO use are not necessarily supported on the red line shown in Figure 1. Additionally, note that DROTRIMM achieves a violation probability that is always lower than or equal to that attained by KNNDRO. Moreover, if we consider values of the robustness parameter smaller than $10^{-1}$ in Figure 2c, many of the dispatch solutions delivered by DROTRIMM, DROW, and KNNDRO lead to a violation probability higher than the threshold $\epsilon = 0.1$, which highlights the value of *distributional robustness* in the OPF problem under uncertainty.

When $N = 30$, DROW exhibits a remarkably robust behavior in the sense that the OPF solutions that this method delivers exhibit comparatively lower variance in terms of both expected cost and reliability. However, like DROTRIMM, DROW also needs a value of its robustness parameter around $10^{-1}$ (or greater) to ensure that almost none of the dispatch solutions it provides violate the reliability threshold $\epsilon = 0.1$. In fact, it is within a neighborhood of that value for the robustness parameter where

the three methods provide reliable OPF solutions with the best cost performance on average. Nevertheless, DROTRIMM is able to identify dispatch solutions that are about 1% cheaper in expectation than those given by DROW and KNNDRO. In this regard, Figure 3 reveals that DROTRIMM achieves this percentage point in expected cost savings through a slightly different generators' dispatch and participation factors. Indeed, see how close to one another the colored bold lines in that figure are for a value of the robustness parameter around $10^{-1}$. Nonetheless, it is also noteworthy that from $10^{-1}$ onward, the OPF solutions given by DROTRIMM start to diverge from those provided by DROW and KNDDRO in a noticeable way.

DROTRIMM's solutions are systematically more cost-efficient on average and variance (and even more reliable) than those of KNNDRO. The reason for this is that KNNDRO does not explicitly protect the dispatch against the potential conditional inference error incurred by the local predictive method it relies on.

When the sample size is high enough, e.g., $N = 2000$, Figure 4c shows that all the methods provide reliable dispatch solutions (i.e., power dispatches that ensure the threshold $\epsilon = 0.1$) for any value of their respective robustness parameter. Moreover, the three methods provide the best OPF solutions in terms of expected cost for low values of this parameter. However, the solutions given by DROW are about 2.3% more expensive, because this method procures far more upward reserve capacity than needed (see Figures 4b). For this, as can be seen in Figure 5, DROW provides a dispatch for the generators and values for their participation factors that are starkly different from those given by the other two methods. On the other hand, the performance of DROTRIMM and KNNDRO tends to be similar in a large-sample regime. Indeed, both are able to disclose dispatches that result in similar system operating costs while guaranteeing the reliability of the system, although DROTRIMM is still comparatively better (see Figures 4c and 4d in the range of the robustness parameter

10

between $10^{-5}$ y $10^{-3}$). This result is hardly surprising because, as the sample data grows in size, the uncertainty intrinsic to the conditional inference diminishes.

Since we have detected in our numerical experiments that DROTRIMM always provide dispatches with a performance similar or superior to those given by KNNDRO, we have not discussed the latter in the case study of Section 5.

## 4. Data for the illustrative example (3-bus system)

This appendix compiles data pertaining to the illustrative example that has been presented and discussed in 3.

| Generator index | Bus index | $c_j^D$ | $c_j^U$ | $g_j^{\min}$ | $g_j^{\max}$ |
|:---:|:---:|:---:|:---:|:---:|:---:|
| 1 | 1 | 6 | 3 | 0 | 120 |
| 2 | 2 | 2 | 5 | 0 | 80 |
| 3 | 3 | 4 | 8 | 0 | 100 |

Table 1. Generators' location, power output limits and reserve capacity costs

| Generator index | Slopes $m_s$ | | | Intercepts $n_s$ | | |
|:---:|:---:|:---:|:---:|:---:|:---:|:---:|
| | piece 1 | piece 2 | piece 3 | piece 1 | piece 2 | piece 3 |
| 1 | 22 | 26 | 30 | 0 | -173 | -493 |
| 2 | 29 | 37 | 45 | 0 | -231 | -658 |
| 3 | 38 | 55 | 71 | 0 | -601 | -1715 |

Table 2. Slopes ($m_s$) and intercepts ($n_s$) of the generators' piecewise linear cost functions

| index | from bus | to bus | X (reactance p.u.) | Cap (MW) |
|:-----:|:--------:|:------:|:------------------:|:--------:|
| 1 | 1 | 2 | 0.13 | 100 |
| 2 | 1 | 3 | 0.13 | 100 |
| 3 | 2 | 3 | 0.13 | 100 |

Table 3. Transmission line parameters



\end{document}